\documentclass[letterpaper, 10pt, conference]{ieeeconf}
\IEEEoverridecommandlockouts 
\usepackage{amsmath,amssymb,url}
\usepackage{graphicx,subfigure,warmread}
\usepackage[all,import]{xy}
\usepackage{color}

\newcommand{\norm}[1]{\ensuremath{\left\| #1 \right\|}}

\newcommand{\bracket}[1]{\ensuremath{\left[ #1 \right]}}
\newcommand{\braces}[1]{\ensuremath{\left\{ #1 \right\}}}

\newcommand{\refeqn}[1]{(\ref{eqn:#1})}
\newcommand{\reffig}[1]{Figure \ref{fig:#1}}
\newcommand{\tr}[1]{\mbox{tr}\ensuremath{\negthickspace\bracket{#1}}}

\newcommand{\SO}{\ensuremath{\mathsf{SO(3)}}}
\newcommand{\T}{\ensuremath{\mathsf{T}}}

\newcommand{\so}{\ensuremath{\mathfrak{so}(3)}}
\newcommand{\SE}{\ensuremath{\mathsf{SE(3)}}}

\renewcommand{\Re}{\ensuremath{\mathbb{R}}}
\newcommand{\Sph}{\ensuremath{\mathsf{S}}}

\newcommand{\D}{\ensuremath{\mathbf{D}}}

\title{\LARGE \bf
Control of Complex Maneuvers for a Quadrotor UAV using Geometric Methods on \SE}

\author{Taeyoung Lee\authorrefmark{1}, Melvin Leok\authorrefmark{2}, and N. Harris McClamroch%
\thanks{Taeyoung Lee, Mechanical and Aerospace Engineering, Florida Institute of Technology, Melbourne, FL 39201 {\tt taeyoung@fit.edu}}%
\thanks{Melvin Leok, Mathematics, University of California at San Diego, La Jolla, CA 92093 {\tt mleok@math.ucsd.edu}}%
\thanks{N. Harris McClamroch, Aerospace Engineering, University of Michigan, Ann Arbor, MI 48109 {\tt
nhm@umich.edu}}%
\thanks{\textsuperscript{\footnotesize\ensuremath{*}}This research has been supported in part by NSF under grants CMMI-1029551.}
\thanks{\textsuperscript{\footnotesize\ensuremath{\dagger}}This research has been supported in part by NSF under grants DMS-0726263, DMS-1001521, DMS-1010687, and CMMI-1029445.}
}

\newcommand{\EditTL}[1]{{\color{red}\protect #1}}
\renewcommand{\EditTL}[1]{{\protect #1}}

\newtheorem{definition}{Definition}

\newtheorem{prop}{Proposition}

\begin{document}
\allowdisplaybreaks
\maketitle \thispagestyle{empty} \pagestyle{empty}

\begin{abstract}
This paper provides new results for control of complex flight maneuvers for a quadrotor unmanned aerial vehicle. The flight maneuvers are defined by a concatenation of flight modes, each of which is achieved by a nonlinear controller that solves an output tracking problem. A mathematical model of the quadrotor UAV rigid body dynamics, defined on the configuration space $\SE$, is introduced as a basis for the analysis. We focus on three output tracking problems, namely (1) outputs given by the vehicle attitude, (2) outputs given by the three position variables for the vehicle center of mass, and (3) output given by the three velocity variables for the vehicle center of mass. A nonlinear tracking controller is developed on the special Euclidean group $\SE$ for each flight mode, and the closed loop is shown to have desirable properties that are almost global in each case. Several numerical examples, including one example in which the quadrotor recovers from being initially upside down and another example that includes switching and transitions between different flight modes, illustrate the versatility and generality of the proposed approach. 
\end{abstract}

\section{INTRODUCTION}

A quadrotor unmanned aerial vehicle (UAV) consists of two pairs of counter-rotating rotors and propellers, located at the vertices of a square frame. It is capable of vertical take-off and landing (VTOL), but it does not require complex mechanical linkages, such as swash plates or teeter hinges, that commonly appear in typical helicopters. Due to its simple mechanical structure, it has been envisaged for various applications such as surveillance or mobile sensor networks as well as for educational purposes. 

Despite the substantial interest in quadrotor UAVs, little attention has been paid to constructing nonlinear control systems that can achieve complex aerobatic maneuvers. Linear control systems such as proportional-derivative controllers or linear quadratic regulators are widely used to enhance the stability properties of an equilibrium~\cite{ValBetPAGNCC06,HofHuaAGNCC07,CasLozICSM05,BouMurAR05,Nic04}. The quadrotor dynamics is modeled as a collection of simplified hybrid dynamic modes, where each mode represents a particular local operating region. But, it is required to do complex reachability analyses to  guarantees the safety and performance of such hybrid system~\cite{GilHofIJRR11}. 

A nonlinear controller is developed for the linearized dynamics of a quadrotor UAV in~\cite{GueHamPIICCA05}. Backstepping and sliding mode techniques are applied in~\cite{BouSiePIICRA05}. Since all of these controllers are based on Euler angles, they exhibit singularities when representing complex rotational maneuvers of a quadrotor UAV, thereby significantly restricting their ability to achieve complex flight maneuvers. 

An attitude control system based on quaternions is applied to a quadrotor UAV~\cite{TayMcGITCSTI06}.
Quaternions do not have singularities, but they have ambiguities in representing an attitude, as the three-sphere, the unit-vectors in $\Re^4$, double-covers the attitude configuration of the special orthogonal group, $\SO$. Therefore, a single physical attitude of a rigid body may yields two different control inputs, which causes inconsistency in the resulting control system. A specific choice between two quaternions generates discontinuity that makes the resulting control system sensitive to noise and disturbances~\cite{MaySanPICDC09}. It is possible to construct continuous controllers, but they may exhibit unwinding behavior, where the controller unnecessarily rotates a rigid body through large angles, even if the initial attitude is close to the desired attitude, thereby breaking Lyapunov stability~\cite{BhaBerSCL00}.

Geometric control, as utilized in this paper, is concerned with the development of control systems for dynamic systems evolving on nonlinear manifolds that cannot be globally identified with Euclidean spaces~\cite{Jur97,BulLew05}. By characterizing geometric properties of nonlinear manifolds intrinsically, geometric control techniques provide unique insights into control theory that cannot be obtained from dynamic models represented using local coordinates. This approach has been applied to fully actuated rigid body dynamics on Lie groups to achieve almost global asymptotic stability~\cite{BulLew05,MaiBerITAC06,CabCunPICDC08,ChaMcCITAC09}.

In this paper, we make use of geometric methods to define and analyze controllers that can achieve complex aerobatic maneuvers for a quadrotor UAV. The dynamics of the quadrotor UAV are expressed globally on the configuration manifold, which is the special Euclidean group $\SE$.  Based on a hybrid control architecture, we construct controllers that can achieve output tracking for outputs that correspond to each of several flight modes, namely an attitude controlled flight mode, a position controlled flight mode, and a velocity controlled flight mode.

The proposed controller exhibits the following unique features: (i) It guarantees almost global tracking features of a quadrotor UAV as the region of attraction almost covers the attitude configuration space $\SO$. Such global stability analysis on the special Euclidean group of a quadrotor UAV is unprecedented. (ii) Hybrid control structures between different tracking mode is robust to switching conditions due to the almost global stability properties. Therefore, aggressive maneuvers of a quadrotor UAV can be achieved in a unified way, without need for complex reachability analyses. (iii) The proposed control system extends the existing geometric controls of the rigid body dynamics into an underactuated rigid body system, where its translation dynamics is coupled to the rotational dynamics in a unique way. (iv) It is coordinate-free. Therefore, it completely avoids singularities, complexities, discontinuities, or ambiguities that arise when using local coordinates or quaternions. 

The paper is organized as follows. We develop a globally defined model for the translational and rotational dynamics of a quadrotor UAV in Section \ref{sec:QDM}. The hybrid control architecture and three flight modes are introduced in Section \ref{sec:GTC}.  Section \ref{sec:ACFM} presents results for the attitude controlled flight mode; Sections \ref{sec:PCFM} and \ref{sec:VCFM} present results for the position controlled flight mode, and the velocity controlled flight mode, respectively.    Several numerical results that demonstrate complex aerobatic maneuvers for a typical quadrotor UAV are presented in Section \ref{sec:NE}.

\section{QUADROTOR DYNAMICS MODEL}\label{sec:QDM}

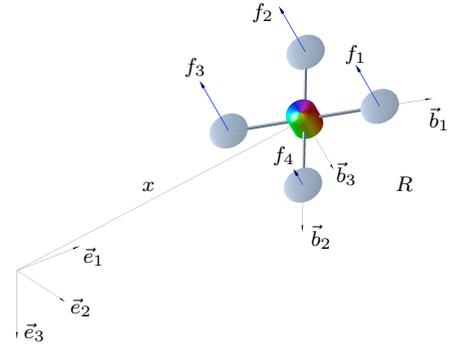
\begin{figure}
\setlength{\unitlength}{0.65\columnwidth}\footnotesize
\centerline{
\begin{picture}(1,0.8)(0,0)
\put(0,0){\includegraphics[width=0.65\columnwidth]{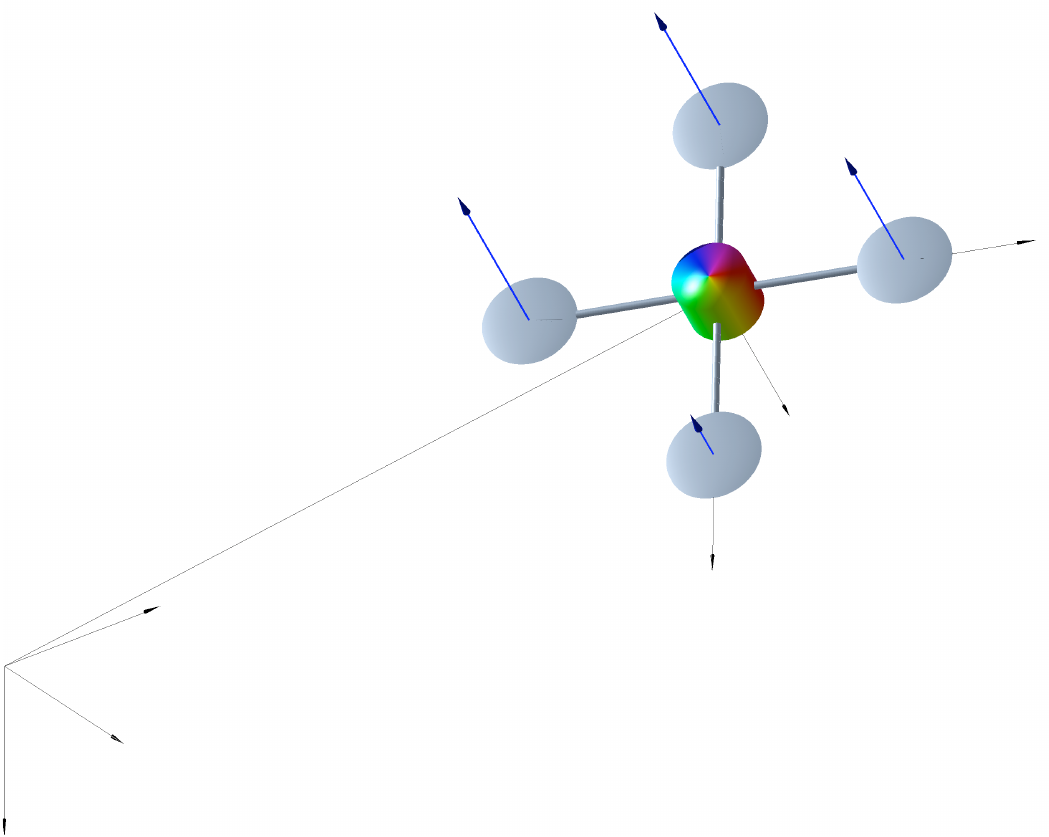}}
\put(0.16,0.18){\shortstack[c]{$\vec e_1$}}
\put(0.13,0.06){\shortstack[c]{$\vec e_2$}}
\put(0.02,0.0){\shortstack[c]{$\vec e_3$}}
\put(0.98,0.5){\shortstack[c]{$\vec b_1$}}
\put(0.70,0.22){\shortstack[c]{$\vec b_2$}}
\put(0.76,0.37){\shortstack[c]{$\vec b_3$}}
\put(0.78,0.66){\shortstack[c]{$f_1$}}
\put(0.56,0.76){\shortstack[c]{$f_2$}}
\put(0.40,0.63){\shortstack[c]{$f_3$}}
\put(0.61,0.42){\shortstack[c]{$f_4$}}
\put(0.30,0.35){\shortstack[c]{$x$}}
\put(0.90,0.35){\shortstack[c]{$R$}}
\end{picture}}
\caption{Quadrotor model}\label{fig:QM}
\end{figure}

Consider a quadrotor UAV model illustrated in \reffig{QM}. This is a system of four identical rotors and propellers located at the vertices of a square, which generate a thrust and torque normal to the plane of this square. We choose an inertial reference frame $\{\vec e_1,\vec e_2,\vec e_3\}$ and a body-fixed frame $\{\vec b_1,\vec b_2,\vec b_3\}$. The origin of the body-fixed frame is located at the center of mass of this vehicle. The first and the second axes of the body-fixed frame, $\vec b_1,\vec b_2$, lie in the plane defined by the centers of the four rotors, as illustrated in \reffig{QM}. The third body-fixed axis $\vec b_3$ is normal to this plane.  Each of the inertial reference frame and the body-fixed reference frame consist of a triad of orthogonal vectors defined according to the right hand rule.    In the subsequent development, these references frames are taken as basis sets and we use vectors in $\Re^3$ to represent physical vectors and we use $3 \times 3$ real matrices to represent linear transformations between the vector spaces defined by these two frames.    Define
\begin{tabular}{lp{5.6cm}}
{$m\in\Re$} & the total mass\\
{$J\in\Re^{3\times 3}$} & the inertia matrix with respect to the body-fixed frame\\
{$R\in\SO$} & the rotation matrix from the body-fixed frame to the inertial  frame\\
{$\Omega\in\Re^3$} & the angular velocity in the body-fixed frame\\
{$x\in\Re^3$} & the position vector of the center of mass in the inertial frame\\
{$v\in\Re^3$} & the velocity vector of the center of mass in the inertial frame\\
{$d\in\Re$} & the distance from the center of mass to the center of each rotor in the $\vec b_1,\vec b_2$ plane\\
{$f_i\in\Re$} & the thrust generated by the $i$-th propeller along the $-\vec b_3$ axis\\
\end{tabular}
\begin{tabular}{lp{5.6cm}}
{$\tau_i\in\Re$} & the torque generated by the $i$-th propeller about the $\vec b_3$ axis\\
{$f\in\Re$} & the total thrust magnitude, i.e., $f=\sum_{i=1}^4 f_i$\\
{$M\in\Re^3$} & the total moment vector in the body-fixed frame\\
\end{tabular}

The configuration of this quadrotor UAV is defined by the location of the center of mass and the attitude with respect to the inertial frame. Therefore, the configuration manifold is the special Euclidean group $\SE$, which is the semidirect product of $\Re^3$ and the special orthogonal group $\SO=\{R\in\Re^{3\times 3}\,|\, R^TR=I,\, \det{R}=1\}$. 

The following conventions are assumed for the rotors and propellers, and the thrust and moment that they exert on the quadrotor UAV. We assume that the thrust of each propeller is directly controlled, and the direction of the thrust of each propeller is normal to the quadrotor plane. The first and third propellers are assumed to generate a thrust along the direction of $-\vec b_3$ when rotating clockwise; the second and fourth propellers are assumed to generate a thrust along the same direction of $-\vec b_3$ when rotating counterclockwise. Thus, the thrust magnitude is $f=\sum_{i=1}^4 f_i$, and it is positive when the total thrust vector acts  along $-\vec b_3$, and it is negative when the total thrust vector acts along $\vec b_3$. By the definition of the rotation matrix $R\in\SO$, 
the total thrust vector is given by $-fRe_3\in\Re^3$ in the inertial frame.
We also assume that the torque generated by each propeller is directly proportional to its thrust. Since it is assumed that the first and the third propellers rotate clockwise and the second and the fourth propellers rotate counterclockwise to generate a positive thrust along the direction of $-\vec b_3$, %
the torque generated by the $i$-th propeller about $\vec b_3$ can be written as $\tau_i=(-1)^{i} c_{\tau f} f_i$  for a fixed constant $c_{\tau f}$.    All of these assumptions are common~\cite{TayMcGITCSTI06,CasLozICSM05}.  The presented control system can readily be extended to include linear rotor dynamics, as studied in~\cite{BouSiePIICRA05}.

Under these assumptions, the moment vector in the body-fixed frame is given by 
\begin{align}
\begin{bmatrix} f \\ M_1 \\ M_2 \\ M_3 $ $\end{bmatrix}
=\begin{bmatrix} 1 & 1& 1& 1\\
0 & -d & 0 & d\\
d & 0 & -d & 0\\
-c_{\tau f} & c_{\tau f} & -c_{\tau f} &c_{\tau f}
\end{bmatrix}
\begin{bmatrix} f_1 \\ f_2\\ f_3 \\ f_4\end{bmatrix}\label{eqn:fM}.
\end{align}
The determinant of the above $4\times4$ matrix is $8 c_{\tau f} d^2$, so it is invertible when $d\neq 0$ and $c_{\tau f}\neq 0$. Therefore, for given thrust magnitude $f$ and given moment vector $M$, the thrust of each propeller $f_1, f_2, f_3, f_4$ can be obtained from \refeqn{fM}. Using this equation, the thrust magnitude $f\in\Re$ and the moment vector $M\in\Re^3$ are viewed as control inputs in this paper.

The equations of motion of the quadrotor UAV can be written as
\begin{gather}
\dot x  = v,\label{eqn:EL1}\\
m \dot v = mge_3 - f R e_3,\label{eqn:EL2}\\
\dot R = R\hat\Omega,\label{eqn:EL3}\\
J\dot \Omega + \Omega\times J\Omega = M,\label{eqn:EL4}
\end{gather}
where the \textit{hat map} $\hat\cdot:\Re^3\rightarrow\so$ is defined by the condition that $\hat x y=x\times y$ for all $x,y\in\Re^3$ (see Appendix \ref{app:hat}). Throughout this paper, $\lambda_m(\cdot)$ and $\lambda_{M}(\cdot)$ denote the minimum eignevalue and the maximum eigenvalue of a matrix, respectively.

\section{GEOMETRIC TRACKING CONTROLS}\label{sec:GTC}

Since the quadrotor UAV has four inputs, it is possible to achieve asymptotic output tracking for at most four quadrotor UAV outputs.    The quadrotor UAV has three translational and three rotational degrees of freedom; it is not possible to achieve asymptotic output tracking of both attitude and position of the quadrotor UAV.    This motivates us to introduce several flight modes.  Each flight mode is associated with a specified set of outputs for which exact tracking of those outputs define that flight mode.   

The three flight modes considered in this paper are:
\begin{itemize}
\item Attitude controlled flight mode: the outputs are the attitude of the quadrotor UAV and the controller for this flight mode achieves asymptotic attitude tracking.  
\item Position controlled flight mode: the outputs are the position vector of the center of mass of the quadrotor UAV and the controller for this flight mode achieves asymptotic position tracking.  
\item Velocity controlled flight mode: the outputs are the velocity vector of the center of mass of the quadrotor UAV and the controller for this flight mode achieves asymptotic velocity tracking.  
\end{itemize}
A complex flight maneuver can be defined by specifying a concatenation of flight modes together with conditions for switching between them; for each flight mode one also specifies the desired or commanded outputs as functions of time.   For example, one might define a complex aerobatic flight maneuver for the quadrotor UAV that consists of a hovering flight segment by specifying a constant position vector, a reorientation segment by specifying the time evolution of the vehicle attitude, and a surveillance flight segment by specifying a time-varying position vector.   The controller in such a case would switch between nonlinear controllers defined for each of the flight modes.   

These types of complex aerobatic maneuvers, involving large angle transitions between flight modes, have not been much studied in the literature.   Such a hybrid flight control architecture has been proposed in~\cite{OiToCDC99,GhToAIAA00,OiToACC02,GilHofIJRR11}, but they are sensitive to switching conditions as the region of attraction for each flight mode is limited, and they required complicated reachability set analyses to guarantee safety and performance. The proposed control system is robust to switching conditions since each flight mode has almost global stability properties, and it is straightforward to design a complex maneuver of a quadrotor UAV.



\section{ATTITUDE CONTROLLED FLIGHT MODE}\label{sec:ACFM}


An arbitrary smooth attitude tracking command $R_d(t) \in \SO$ is given as a function of time. The corresponding angular velocity command is obtained by the attitude kinematics equation, $\hat\Omega_d = R_d^T \dot R_d$. We first define errors associated with the attitude dynamics of the quadrotor UAV.   The attitude and angular velocity tracking error should be carefully chosen as they evolve on the tangent bundle of $\SO$.   First, define the real-valued error function on $\SO \times \SO$: 
\begin{align}                                            
\Psi(R,R_d) = \frac{1}{2}\tr{I-R_d^T R}.\label{eqn:Psi}
\end{align}
This function is locally positive-definite about $R=R_d$ within the region where the rotation angle between $R$ and  $R_d$ is less than $180^\circ$~\cite{BulLew05}.  For a given $R_d$, this set can be represented by the sublevel set $L_{2}=\{R\in\SO\,|\,\Psi(R,R_d) < 2\}$, which almost covers $\SO$. 

The variation of a rotation matrix can be expressed as $\delta R = R\hat\eta$ for $\eta\in\Re^3$, so that the derivative of the error function is given by
\begin{align}
\D_R \Psi(R,R_d)\cdot R\hat\eta & = -\frac{1}{2}\tr{R_d^T R\hat\eta}=e_R\cdot\eta
\end{align}
where the attitude tracking error $e_R\in\Re^3 $ is chosen as
\begin{align}
e_R = \frac{1}{2} (R_d^TR - R^T R_d)^\vee.
\end{align}
The \textit{vee map} $^\vee:\so\rightarrow\Re^3$ is the inverse of the hat map. We used a property of the hat map given by equation  \refeqn{hat1} in Appendix A. 

The tangent vectors $\dot R\in\T_R\SO$ and $\dot R_d\in\T_{R_d}\SO$ cannot be directly compared since they lie in different tangent spaces. We transform $\dot R_d$ into a vector in $\T_{R}\SO$, and we compare it with $\dot R$ as follows:
\begin{align*}
\dot R - \dot R_d (R_d^T R)  
& = R(\hat\Omega - R^TR_d\hat\Omega_d R_d^T R) = R\hat e_\Omega,
\end{align*}
where the tracking error for the angular velocity $e_\Omega\in\Re^3$ is defined as follows:
\begin{align}
e_\Omega = \Omega - R^T R_d \Omega_d.\label{eqn:eW}
\end{align}
We show that $e_\Omega$ is the angular velocity vector of the relative rotation matrix $R_d^T R$, represented in the body-fixed frame, since
\begin{align}
\frac{d}{dt} (R_d^T R) 
= (R_d^T R)\, \hat e_\Omega.\label{eqn:Redot}
\end{align}

We now introduce a nonlinear controller for the attitude controlled flight mode, described by an expression for the moment vector:
\begin{align}
M & = -k_R e_R -k_\Omega e_\Omega +\Omega\times J\Omega\nonumber\\
&\qquad -J(\hat\Omega R^T R_d \Omega_d - R^T R_d\dot\Omega_d),\label{eqn:aM}
\end{align}
where $k_R,k_\Omega$ are positive constants. 

In this attitude controlled mode, it is possible to ignore the translational motion of the quadrotor UAV; consequently the reduced model for the attitude dynamics are given by equations \refeqn{EL3}, \refeqn{EL4}, using the controller expression \refeqn{aM}.  We now state the result that $(e_R, e_\Omega) = (0,0)$ is an exponentially stable equilibrium of the reduced closed loop dynamics.  

\begin{prop}{(Exponential Stability of Attitude Controlled Flight Mode)}\label{prop:Att}
Consider the control moment $M$ defined in \refeqn{aM} for any positive constants $k_R,k_\Omega$. 
Suppose that the initial conditions satisfy
\begin{gather}
\Psi(R(0),R_d(0))<2,\label{eqn:eRb0}\\
\norm{e_\Omega(0)}^2 < \frac{2}{\lambda_{M}(J)} k_R (2-\Psi(R(0),R_d(0))).\label{eqn:eWb}
\end{gather}
Then, the zero equilibrium of the closed loop tracking error $(e_R,e_\Omega)=(0,0)$ is exponentially stable.
 Furthermore, there exist constants $\alpha_2,\beta_2 >0$ such that
\begin{align}
\Psi(R(t),R_d(t)) \leq \min\braces{ 2, \alpha_2 e^{-\beta_2 t}}.\label{eqn:Psib}
\end{align}
\end{prop}

\begin{proof}
See Appendix \ref{sec:pfAtt}.
\end{proof}

In this proposition, equations \refeqn{eRb0}, \refeqn{eWb}  describe a region of attraction for the reduced closed loop dynamics.   An estimate of the domain of attraction is obtained for which the quadrotor attitude lies in the sublevel set $L_2=\{R\in\SO\,|\,\Psi(R,R_d) < 2\}$ for a given $R_d$.   This requires that the initial attitude error should be less than $180^\circ$, in terms of the rotation angle about the eigenaxis between $R$ and $R_d$. Therefore, in Proposition \ref{prop:Att}, exponential stability is guaranteed for almost all initial attitude errors. More explicitly, the attitudes that lie outside of the region of attraction are of the form $\exp(\pi \hat s)R_d$ for some $s\in\Sph^2$. Since they comprise a two-dimensional manifold in the three-dimensional $\SO$, we claim that the presented controller exhibits \textit{almost global} properties in $\SO$.
%
%
It should be noted that topological obstructions prevent one from constructing a smooth controller on $\SO$ that has an equilibrium solution that is global asymptotically stable~\cite{KodPICDC98}. The size of the region of attraction 
can be increased by choosing a larger controller gain $k_R$ in \refeqn{eWb}. 

Asymptotic tracking of the quadrotor attitude does not require specification of the thrust magnitude.   As an auxiliary problem, the thrust magnitude can be chosen in many different ways to achieve an additional translational motion objective. 

As an example of a specific selection approach, we assume that the objective is to asymptotically track a quadrotor altitude command.    It is straightforward to obtain the following corollary of Proposition 1.

\begin{prop}{(Exponential Stability of Attitude Controlled Flight Mode with Altitude Tracking)}\label{prop:Alt}
Consider the control moment vector $M$ defined in \refeqn{aM} satisfying the assumptions of Proposition 1.  In addition, the thrust magnitude is given by
\begin{gather}
f  =  \frac{k_x ( x_3 - x_{3_d}) + k_v (\dot x_3 - \dot x_{3_d}) + mg - m\ddot x_{3_d}}{e_3\cdot R e_3},\label{eqn:af}
\end{gather}
where $k_x, k_v$ are positive constants, $x_{3_d}(t)$  is the quadrotor altitude command, and we assume that 
\begin{gather}
e_3\cdot R e_3 \not = 0.
\end{gather}
The conclusions of Proposition 1 hold and in addition the quadrotor altitude $x_3(t)$ asymptotically tracks the altitude command $x_{3_d}(t)$.   
\end{prop}
\begin{proof}
See Appendix \ref{sec:pfAlt}.
\end{proof}


Since the translational motion of the quadrotor UAV can only be partially controlled; this flight mode is most suitable for short time periods where an attitude maneuver is to be completed.   The translational equations of motion of the quadrotor UAV, during an attitude flight mode, are given by equations \refeqn{EL1}, \refeqn{EL2}, and whatever thrust magnitude controller, e.g., equation \refeqn{af}, is selected. 

\section{POSITION CONTROLLED FLIGHT MODE}\label{sec:PCFM}

We now introduce a nonlinear controller for the position controlled flight mode.  We show that this controller achieves almost global asymptotic position tracking, that is the output position vector of the quadrotor UAV asymptotically tracks the commanded position.     This flight mode requires analysis of the coupled translational and rotational equations of motion; hence, we make use of the notation and analysis in the prior section to describe the properties of the closed loop system in this flight mode.  

An arbitrary smooth position tracking command $x_d(t) \in \Re^3$ is chosen.    The position tracking errors for the position and the velocity are given by:
\begin{align}
e_x & = x - x_d,\\
e_v & = v - \dot x_d.
\end{align}

The nonlinear controller for the position controlled flight mode, described by control expressions for the  thrust magnitude and the  moment vector, are:
\begin{align}
f & = ( k_x e_x + k_v e_v + mg e_3-m\ddot x_d)\cdot Re_3,\label{eqn:f}\\
M & = -k_R e_R -k_\Omega e_\Omega +\Omega\times J\Omega\nonumber\\
&\quad -J(\hat\Omega R^T R_c \Omega_c - R^T R_c\dot\Omega_c),\label{eqn:M}
\end{align}
where $k_x,k_v,k_R,k_\Omega$ are positive constants.  Following the prior definition of the attitude error and the angular velocity error
\begin{align}
e_R = \frac{1}{2} (R_c^TR - R^T R_c)^\vee, \quad 
e_\Omega = \Omega - R^T R_c \Omega_c\label{eqn:eWc},
\end{align}
and the computed attitude $R_c(t)\in\SO$ and computed angular velocity $\Omega_c \in \mathbb{R}^3$ are given by
\begin{align}
R_c=[ b_{1_c};\, b_{3_c}\times b_{1_c};\, b_{3_c}],\quad \hat\Omega_c = R_c^T \dot R_c\label{eqn:RdWc},
\end{align}
where $ b_{3_c} \in \Sph^2$ is defined by
\begin{align}
 b_{3_c} = -\frac{-k_x e_x - k_v e_v - mg e_3 +m\ddot x_d}{\norm{-k_x e_x - k_v e_v - mg e_3 + m\ddot x_d}},\label{eqn:Rd3}
\end{align}
and $ b_{1c} \in \Sph^2$ is selected to be orthogonal to $ b_{3c}$, thereby guaranteeing that $R_c \in \SO$. We assume that 
\begin{align}
\norm{-k_x e_x - k_v e_v - mg e_3 + m\ddot x_d} \neq 0,\label{eqn:A1}
\end{align}
and the commanded acceleration is uniformly bounded such that
\begin{align}
\|-mge_3+m\ddot x_d\| < B\label{eqn:B}
\end{align}
for a given positive constant $B$. 

The thrust magnitude controller and the moment vector controller is feedback dependent on the position and translational velocity and they depend on the commanded position, translational velocity and translational acceleration.    The control moment vector has a form that is similar to that for the attitude controlled flight mode.   However, the \textit{attitude error} and \textit{angular velocity error} are defined with respect to a computed attitude, angular velocity and angular acceleration, that are constructed according to the indicated procedure.  

The nonlinear controller given by equations \refeqn{f}, \refeqn{M} can be given a backstepping interpretation.   The computed attitude $R_c$ given in equation \refeqn{RdWc} is selected so that the thrust axis $-b_3$ of the quadrotor UAV tracks the computed direction given by $-b_{3_c}$ in \refeqn{Rd3}, which is a direction of the thrust vector that achieves position tracking.   The moment expression \refeqn{M} causes the attitude of the quadrotor UAV to asymptotically track $R_c$ and the thrust magnitude expression \refeqn{f} achieves asymptotic position tracking.  

\begin{figure}
\centerline{
\setlength{\unitlength}{2.1em}\centering\footnotesize
\begin{picture}(11,3.2)(0.0,-3.2)
\put(1.2,-1.5){\framebox(2,1.0)[c]
{\shortstack[c]{Force\\controller}}}
\put(4.0,-2.2){\framebox(2,1.0)[c]
{\shortstack[c]{Moment\\controller}}}
\put(0,-0.95){\vector(1,0){1.2}}
\put(3.2,-1.35){\vector(1,0){0.8}}
\put(3.2,-0.65){\vector(1,0){4.1}}
\put(0,-1.8){\vector(1,0){4.0}}
\put(6,-1.7){\vector(1,0){1.3}}
\put(7.3,-2.25){\framebox(2.5,1.8)[c]
{\shortstack[c]{Quadrotor\\Dynamics}}}
\put(9.8,-1.35){\vector(1,0){1.2}}
\put(6.6,-1.0){$f$}
\put(6.5,-2.05){$M$}
\put(3.30,-1.15){$b_{3_c}$}
\put(0.25,-0.80){$x_d$}
\put(0.0,-1.65){($b_{1_d}$)}
\put(5.5,-3.2){$x,v,R,\Omega$}
\put(10.3,-1.35){\line(0,-1){1.5}}
\put(10.3,-2.85){\line(-1,0){8.1}}
\put(2.2,-2.85){\vector(0,1){1.35}}
\put(2.2,-2.10){\vector(1,0){1.8}}
\put(2.2,-2.10){\circle*{0.1}}
\put(10.3,-1.35){\circle*{0.1}}
\put(0.9,-2.4){\dashbox{0.08}(5.4,2.1)[c]{}}
\put(2.8,-2.7){Controller}
\end{picture}
}
\caption{Controller structure for position controlled flight mode}\label{fig:CS}
\end{figure}
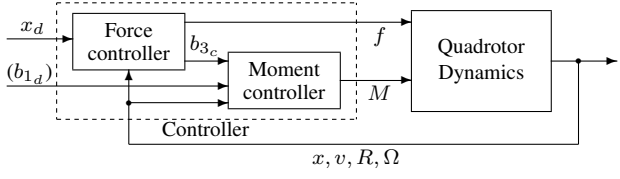

The closed loop system for this position controlled flight mode is illustrated in \reffig{CS}.   The corresponding closed loop control system is described by equations \refeqn{EL1}, \refeqn{EL2}, \refeqn{EL3}, \refeqn{EL4}, using the controller expressions \refeqn{f} and \refeqn{M}. 

We now state the result that $(e_x, e_v, e_R, e_\Omega) = (0,0,0,0)$ is an exponentially stable equilibrium of the closed loop dynamics.

\begin{prop}{(Exponential Stability of Position Controlled Flight Mode)}\label{prop:Pos}
Consider the thrust magnitude $f$ and moment vector $M$ defined by equations \refeqn{f}, \refeqn{M}.   Suppose that the initial conditions satisfy
\begin{gather}
\Psi(R(0),R_c(0)) < 1,\label{eqn:Psi0}\\
\|e_x(0)\| < e_{x_{\max}},\label{eqn:ex0}
\end{gather}
for a fixed constant $e_{x_{\max}}$. Define $W_1,W_{12},W_2\in\Re^{2\times 2}$ to be
\begin{align}
W_1 &= \begin{bmatrix} \frac{c_1k_x}{m}(1-\alpha) & -\frac{c_1k_v}{2m}(1+\alpha)\\
-\frac{c_1k_v}{2m}(1+\alpha) & k_v(1-\alpha)-c_1\end{bmatrix},\label{eqn:W1}\\
W_{12}&=\begin{bmatrix}
\frac{c_1}{m}B & 0 \\ B+k_x e_{x_{\max}} & 0\end{bmatrix},\label{eqn:W12}\\
W_2 &= \begin{bmatrix} \frac{c_2k_R}{\lambda_{M}(J)} & -\frac{c_2k_\Omega}{2\lambda_{m}(J)} \\ 
-\frac{c_2k_\Omega}{2\lambda_{m}(J)} & k_\Omega-c_2 \end{bmatrix},\label{eqn:W20}
\end{align}
where $\Psi(R(0),R_c(0))<\psi_1 < 1$, and $\alpha=\sqrt{\psi_1(2-\psi_1)}$. For positive constants $k_x,k_v$, we choose positive constants $c_1,c_2,k_R,k_\Omega$ such that
\begin{gather}
c_1 < \min\braces{k_v(1-\alpha),\;\frac{4mk_xk_v(1-\alpha)^2}{k_v^2(1+\alpha)^2+4m k_x(1-\alpha)},\; \sqrt{k_xm} },\label{eqn:c1b}\\
c_2 < \min\bigg\{ k_\Omega, \frac{4k_\Omega k_R\lambda_{m}(J)^2}{k_\Omega^2\lambda_{M}(J)+4k_R\lambda_{m}(J)^2}, \sqrt{k_R\lambda_{m}(J)}\bigg\},\label{eqn:c2b}\\
\lambda_{m}(W_2) > \frac{4\|W_{12}\|^2}{\lambda_{m}(W_1)}.\label{eqn:kRkWb}
\end{gather}
Then, the zero equilibrium of the closed loop tracking errors $(e_x,e_v,e_R,e_\Omega)=(0,0,0,0)$ is exponentially stable. A region of attraction is characterized by \refeqn{Psi0}, \refeqn{ex0}, and 
\begin{gather}
\norm{e_\Omega(0)}^2 < \frac{2}{\lambda_{M}(J)} k_R (\psi_1-\Psi(R(0),R_c(0))),\label{eqn:eWb2}\\
\lambda_M(M_{12}) \|z_1(0)\|^2 + \lambda_M(M_{22}')\|z_2(0)\|^2 < \frac{1}{2}k_x e_{x_{\max}}^2,\label{eqn:RAz}
\end{gather}
where $z_1=[\|e_x\|,\,\|e_v\|]^T,\,z_2=[\|e_R\|,\,\|e_\Omega\|]^T\in\Re^2$ and 
\begin{align*}
M_{12} = \frac{1}{2}\begin{bmatrix} k_x & c_1 \\ c_1 & m\end{bmatrix},\quad
M'_{22} = \frac{1}{2}\begin{bmatrix} \frac{2k_R}{2-\psi_1} & c_2 \\ c_2 & \lambda_{M}(J)\end{bmatrix}.
\end{align*}
\end{prop}

\begin{proof}
See Appendix \ref{sec:pfPos}.
\end{proof}



Proposition \ref{prop:Pos} requires that the initial \textit{attitude error} is less than $90^\circ$ to achieve  exponential stability for this flight mode. Suppose that this is not satisfied, i.e. $1\leq\Psi(R(0),R_c(0))<2$.    We can apply Proposition \ref{prop:Att}, which states that the attitude error function $\Psi$ exponentially decreases, and therefore, it enters the region of attraction of Proposition \ref{prop:Pos} in a finite time. Therefore, by combining the results of Proposition \ref{prop:Att} and \ref{prop:Pos}, we can show almost global exponential attractiveness when $\Psi(R(0),R_c(0))<2$.

\begin{definition}{(Exponential Attractiveness~\cite{Qu98})} An equilibrium point $z=0$ of a dynamic systems is \textit{exponentially attractive} if, for some $\delta >0$, there exists a constant $\alpha(\delta) >0$ and $\beta >0$ such that $\norm{z(0)} < \delta$ implies  $\|z(t)\| \leq \alpha(\delta) e^{-\beta t}$ for all $t>0$.
\end{definition}

This should be distinguished from the stronger notion of exponential stability, in which the above bound is replaced by $\|z(t)\| \leq \alpha(\delta)\norm{z(0)} e^{-\beta t}$.

\begin{prop}{(Almost Global Exponential Attractiveness of the Position Controlled Flight Mode)}\label{prop:Pos2}
Consider the thrust magnitude $f$ and moment vector $M$ defined in expressions \refeqn{f}, \refeqn{M}.  Suppose that the initial conditions satisfy
\begin{gather}
1\leq \Psi(R(0),R_c(0))<2\label{eqn:eRb3},\\ \norm{e_\Omega(0)}^2 < \frac{2}{\lambda_{M}(J)} k_R (2-\Psi(R(0),R_c(0))).\label{eqn:eWb3}
\end{gather}
Then, the zero equilibrium of the closed loop tracking errors $(e_x,e_v,e_R,e_\Omega)=(0,0,0,0)$ is exponentially attractive.
\end{prop}

\begin{proof}
See Appendix \ref{sec:pfPos2}.
\end{proof}

In Proposition \ref{prop:Pos2}, exponential attractiveness is guaranteed for almost all initial attitude errors. Since the attitudes that lie outside of the region of attraction comprise a two-dimensional manifold in the three-dimensional $\SO$, as discussed in Section \ref{sec:ACFM}, we claim that the presented controller exhibits \textit{almost global} properties in $\SO$.



As described above, the construction of the orthogonal matrix $R_c$ involves having its third column $b_{3_c}$ specified by a normalized feedback function, and its first column $b_{1_c}$ is chosen to be orthogonal to the third column. The unit vector $b_{1_c}$ can be arbitrarily chosen in the plane normal to $b_{3_c}$, which corresponds to a one-dimensional degree of choice. This reflects the fact that 
the quadrotor UAV has four control inputs that are used to track a three-dimensional position command.

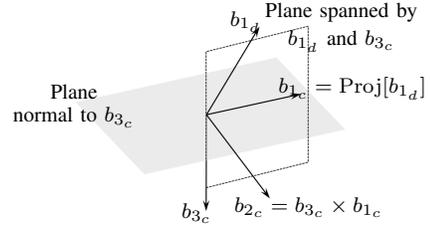
\begin{figure}
\centerline{
\renewcommand{\xyWARMinclude}[1]{\includegraphics[width=0.42\columnwidth]{#1}}
{\footnotesize\selectfont
$$\begin{xy}
\xyWARMprocessEPS{b1d}{pdf}
\xyMarkedImport{}
\xyMarkedMathPoints{1-15}
\end{xy}
$$}}
\caption{Convergence property of the first body-fixed axis: $b_{3_c}$ is determined by \refeqn{Rd3}. 
We choose an arbitrary $b_{1_d}$ that is not parallel to $b_{3_c}$, and project it on to the plane normal to $b_{3_c}$ to obtain $b_{1_c}$. This guarantees that the first body-fixed axis asymptotically lies in the plane spanned by $b_{1_d}$ and $b_{3_c}$, which converges to the direction of $ge_3-\ddot x_d$ as $t\rightarrow\infty$.}\label{fig:b1d}
\vspace*{-0.4cm}
\end{figure}

By choosing $b_{1_c}$ properly, we constrain the asymptotic direction of the first body-fixed axis. 
Here, we propose to specify the \textit{projection} of the first body-fixed axis onto the plane normal to $b_{3_c}$. 
In particular, we choose a desired direction $b_{1_d}\in\Sph^2$, that is not parallel to $b_{3_c}$, and $b_{1_c}$ is selected as
$b_{1_c}=\mathrm{Proj}[b_{1_d}]$, where $\mathrm{Proj}[\cdot]$ denotes the normalized projection onto the plane perpendicular to $b_{3_c}$. In this case, the first body-fixed axis does not converge to $b_{1_d}$, but it converges to the projection of $b_{1_d}$, i.e. $b_1\rightarrow b_{1_c}=\mathrm{Proj}[b_{1_d}]$ as $t\rightarrow\infty$. In other words, the first body-fixed axis converges to the intersection of 
the plane normal to $b_{3_c}$ and the plane spanned by $b_{3_c}$ and $b_{1_d}$ (see \reffig{b1d}). From \refeqn{Rd3}, we observe that $b_{3_c}$ asymptotically converges to the direction $ge_3-\ddot x_d$. In short, the additional input is used to guarantee that the first body-fixed axis asymptotically lies in the plane spanned by $b_{1_d}$ and $ge_3-\ddot x_d$.

Suppose that $\ddot x_d =0$, then the third body-fixed axis converges to the gravity direction $e_3$. In this case, we can choose $b_{1_d}$ arbitrarily in the horizontal plane, and it follows that $b_{1_c}\rightarrow\mathrm{Proj}[b_{1_d}]=b_{1_d}$ as $t\rightarrow\infty$. Therefore, the first body-fixed axis $b_1$ asymptotically converges to $b_{1_d}$, which can be used to specify the heading direction of the quadrotor UAV in the horizontal plane. 

\begin{prop}{(Almost Global Exponential Attractiveness of Position Controlled Flight Mode with Specified Asymptotic Direction of First Body-Fixed Axis)}\label{prop:34C}
Consider the moment vector $M$ defined in \refeqn{M} and the thrust magnitude $f$ defined in \refeqn{f} satisfying the assumptions of Propositions \ref{prop:Pos} and \ref{prop:Pos2}.  

In addition, the first column of $R_c$, namely $b_{1_c}$ is constructed as follows. We choose $b_{1_d}(t)\in\Sph^2$, and we assume that it is not parallel to $b_{3_c}$. The unit vector $b_{1_c}$ is constructed by projecting $b_{1_d}$ onto the plane normal to $b_{3_c}$, and normalizing it:
\begin{align}
b_{1_c} = -\frac{1}{\|b_{3_c}\times b_{1_d}\|} (b_{3_c}\times (b_{3_c}\times b_{1_d})).\label{eqn:b1c}
\end{align}
Then, the conclusions of Propositions \ref{prop:Pos} and \ref{prop:Pos2} hold, and the first body-fixed axis asymptotically lies in the plane spanned by $b_{1_d}$ and $g e_3-\ddot x_d$.

In the special case where $\ddot x_d=0$, we can choose $b_{1_d}$ in the horizontal plane. Then, the first body-fixed axis asymptotically converges to $b_{1_d}$.
\end{prop}

Expressions for $\Omega_c$ and $\dot\Omega_c$ that appear in Proposition \ref{prop:34C} are summarized in~\cite{LeeLeo10}.
These additional properties of the closed loop can be interpreted as characterizing the asymptotic direction of the first body-fixed axis and the asymptotic direction of the third body-fixed axis as it depends on the commanded vehicle acceleration.   These physical properties may be of importance in some flight maneuvers. 


\section{VELOCITY CONTROLLED FLIGHT MODE}\label{sec:VCFM}

We now introduce a nonlinear controller for the velocity controlled flight mode. An arbitrary velocity tracking command $t\rightarrow v_d(t) \in \Re^3$ is given. The velocity tracking error is given by:
\begin{align}
e_v & = v - v_d.
\end{align}

The nonlinear controller for the velocity controlled flight mode is given by
\begin{align}
f & = ( k_v e_v + mg e_3-m\dot v_d)\cdot Re_3,\label{eqn:fV}\\
M & = -k_R e_R -k_\Omega e_\Omega +\Omega\times J\Omega\nonumber\\
&\quad -J(\hat\Omega R^T R_c \Omega_c - R^T R_c\dot\Omega_c),\label{eqn:MV}
\end{align}
where $k_v,k_R,k_\Omega$ are positive constants, and following the prior definition of the attitude error and the angular velocity error
\begin{align}
e_R = \frac{1}{2} (R_c^TR - R^T R_c)^\vee, \quad 
e_\Omega = \Omega - R^T R_c \Omega_c\label{eqn:eWcV},
\end{align}
and $R_c(t) \in \SO$ and $\Omega_c \in R^3$ are constructed as:
\begin{align}
R_c=[ b_{1_c};\, b_{3_c}\times  b_{1_c};\,  b_{3_c}],\quad \hat\Omega_c = R_c^T \dot R_c\label{eqn:RdWcV},
\end{align}
where $ b_{3_c} \in \Sph^2$ is defined by
\begin{align}
 b_{3_c} = -\frac{ - k_v e_v - mg e_3 +m\dot v_d}{\norm{ - k_v e_v - mg e_3 + m\dot v_d}}.\label{eqn:Rd3V}
\end{align}
and $ b_{1c} \in \Sph^2$ is selected to be orthogonal to $ b_{3c}$, thereby guaranteeing that $R_c \in \SO$. We assume that 
\begin{gather}
\norm{- k_v e_v - mg e_3 + m\ddot x_d} \neq 0,\label{eqn:A1V}\\
\|- mg e_3 + m\ddot x_d\| < B\label{eqn:BV}
\end{gather}
for a given positive constant $B$. 

The overall controller structure and the corresponding stability properties are similar to the position controlled flight mode.  More explicitly, the closed loop dynamics have the property that  $(e_v, e_R, e_\Omega) = (0,0,0)$ is an equilibrium that is exponentially stable for any initial condition satisfying $\Psi(R(0),R_c(0)) < 1$, and it is exponentially attractive  for any initial condition satisfying $\Psi(R(0),R_c(0)) < 2$. Due to page limitations, the explicit statements of propositions and proofs for the velocity controlled flight mode are relegated to~\cite{LeeLeo10}.

\section{NUMERICAL RESULTS ILLUSTRATING COMPLEX FLIGHT MANEUVERS}\label{sec:NE}

Numerical results are presented to demonstrate the prior approach for performing complex flight maneuvers for a typical quadrotor UAV.     The parameters are chosen to match a quadrotor UAV described in~\cite{PouMahACRA06}.
\begin{gather*}
J=[0.0820,0.0845,0.1377]\,\mathrm{kg-m^2},\quad m=4.34\,\mathrm{kg}\\
 d=0.315\,\mathrm{m},\quad c_{\tau f}=8.004\times 10^{-3}\,\mathrm{m}.
\end{gather*}
The controller parameters are chosen as follows:
\begin{align*}
k_x=16m,\quad k_v=5.6m,\quad k_R=8.81,\quad k_\Omega = 2.54.
\end{align*}


We consider two complex flight maneuvers.  The first case corresponds to the position controlled mode; the results in Proposition \ref{prop:Pos2} are referenced. The second case involves transitions between all of the three flight modes. 


\paragraph*{Case (I): Position Controlled Flight Mode}
Consider a hovering maneuver for which the quadrotor UAV recovers from being initially upside down.   The desired tracking commands are as follows.
\begin{gather*}
x_d(t) = [0, 0, 0],\quad b_{1_d} (t) = [1,0 ,0].
\end{gather*}
and it is desired to maintain the quadrotor UAV at a constant altitude.   
Initial conditions are chosen as
\begin{gather*}
x(0)=[0,0,0],\quad v(0)=[0,0,0],\\
R(0)=\begin{bmatrix}
            1&            0&            0\\
            0&  -0.9995 &  -0.0314\\
            0&   0.0314 & -0.9995\end{bmatrix},\quad
\Omega(0)=[0,0,0].                       
\end{gather*}
This initial condition corresponds to an upside down quadrotor UAV.  

The preferred direction of the total thrust vector in the controlled system is $-b_3$. But initially, it is given by $-b_3(0)=-R(0)e_3 = [0,0.0314,0.9995]$, which is almost opposite to the thrust direction $[0,0,-1]$ required for the given hovering command. This yields a large initial attitude error, namely $178^\circ$ in terms of the rotation angle about the eigen-axis between $R_c(0)$ and $R(0)$, and the corresponding the initial \textit{attitude error} is $\Psi(0)=1.995$. 

\begin{figure}
\centerline{
	\subfigure[Attitude error function $\Psi$]{
		\includegraphics[width=0.48\columnwidth]{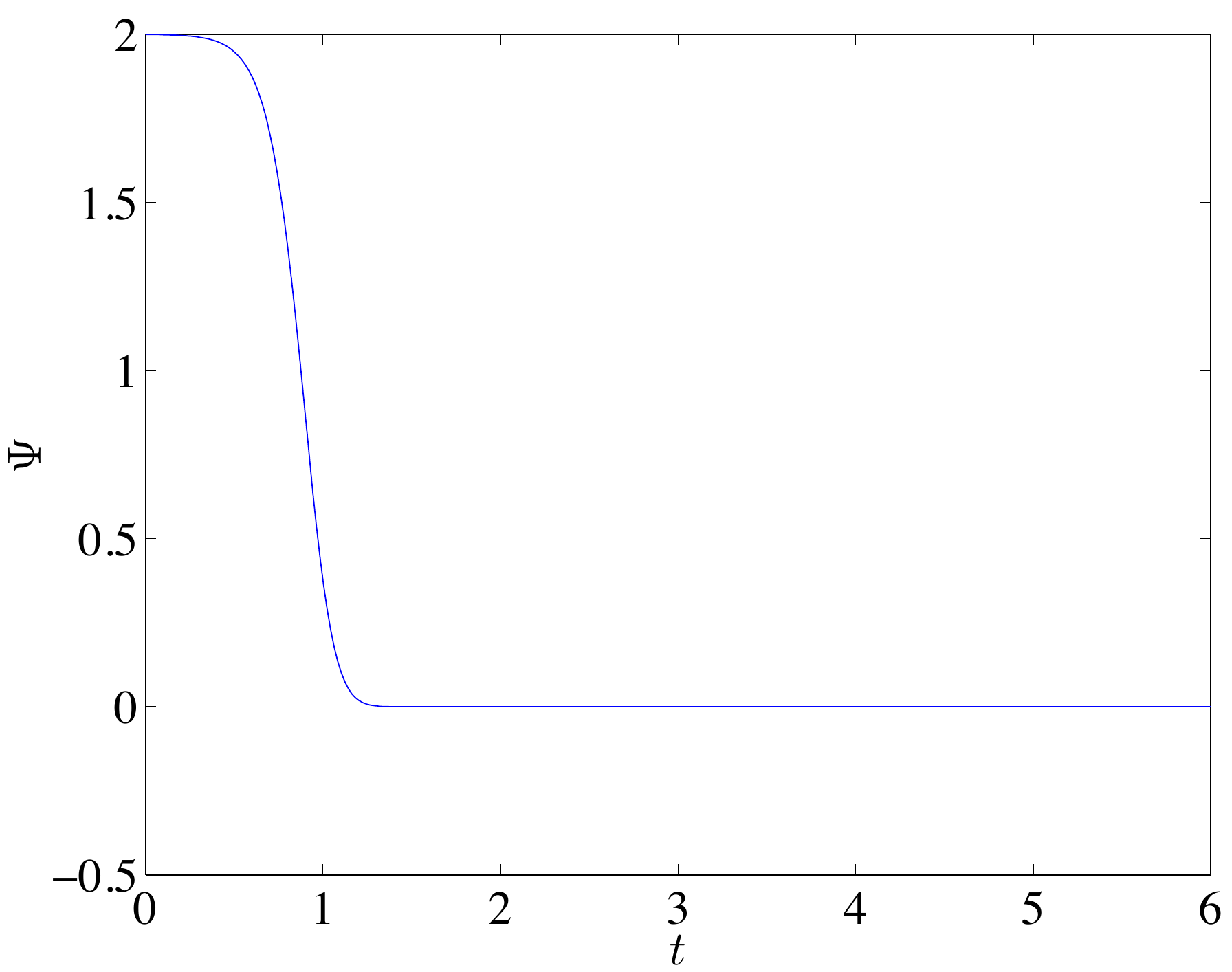}\label{fig:IIPsi}}
		\hfill
	\subfigure[Position ($x$:solid, $x_d$:dotted, ($\mathrm{m}$))]{
		\includegraphics[width=0.48\columnwidth]{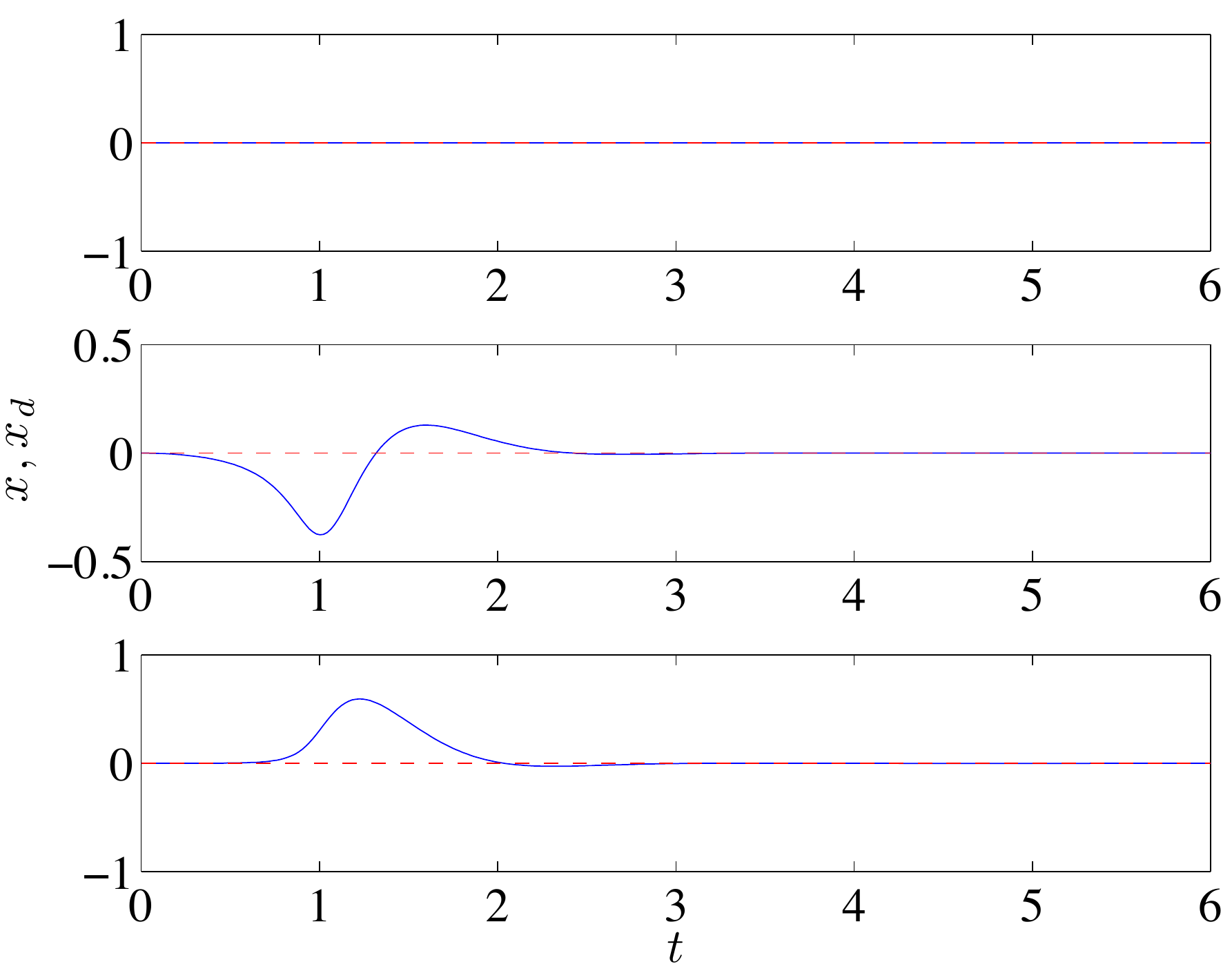}\label{fig:IIx}}
}
\centerline{
	\subfigure[Angular velocity ($\mathrm{rad/sec}$)]{
		\includegraphics[width=0.48\columnwidth]{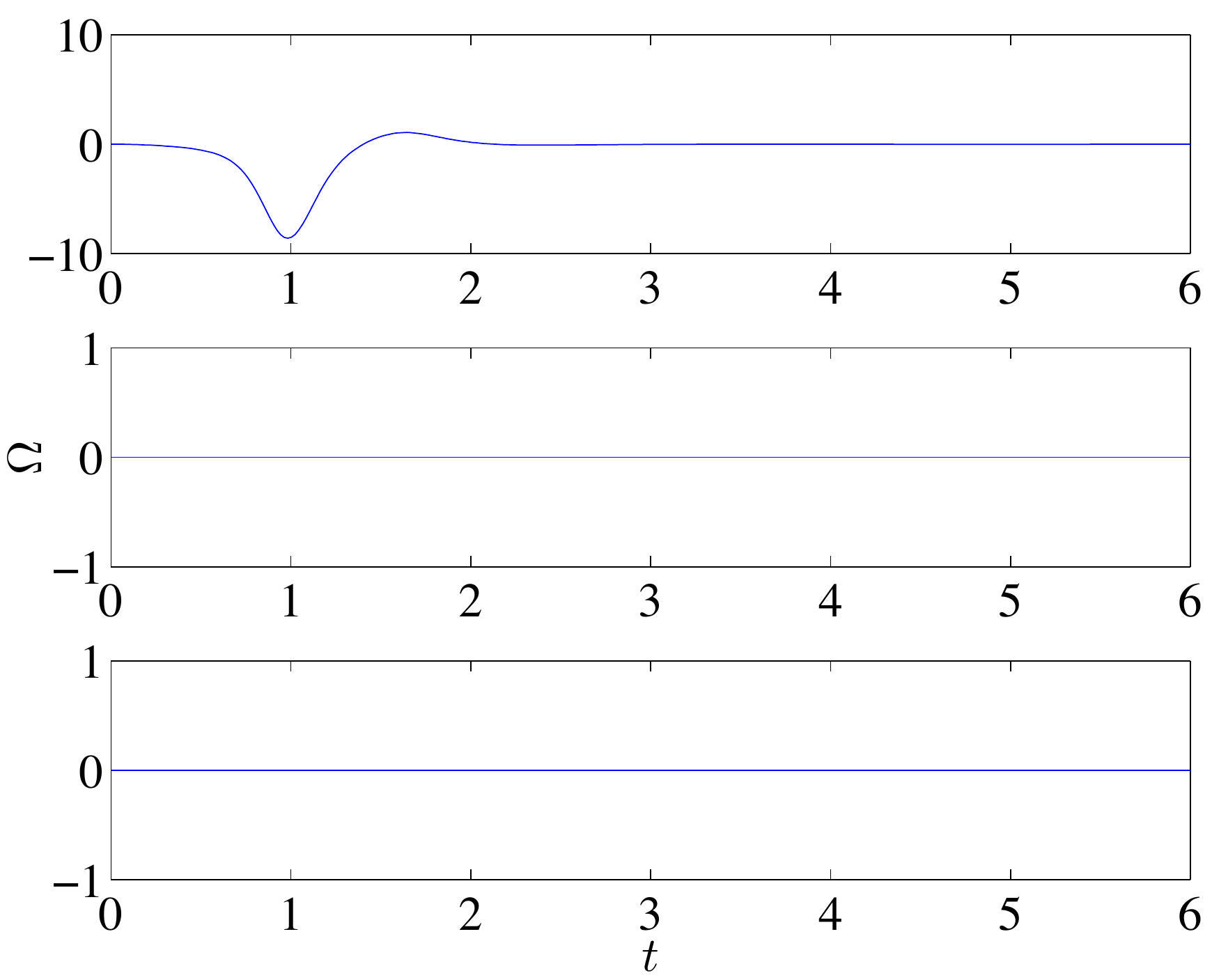}\label{fig:IIW}}
		\hfill
	\subfigure[Thrust of each rotor ($\mathrm{N}$)]{
		\includegraphics[width=0.50\columnwidth]{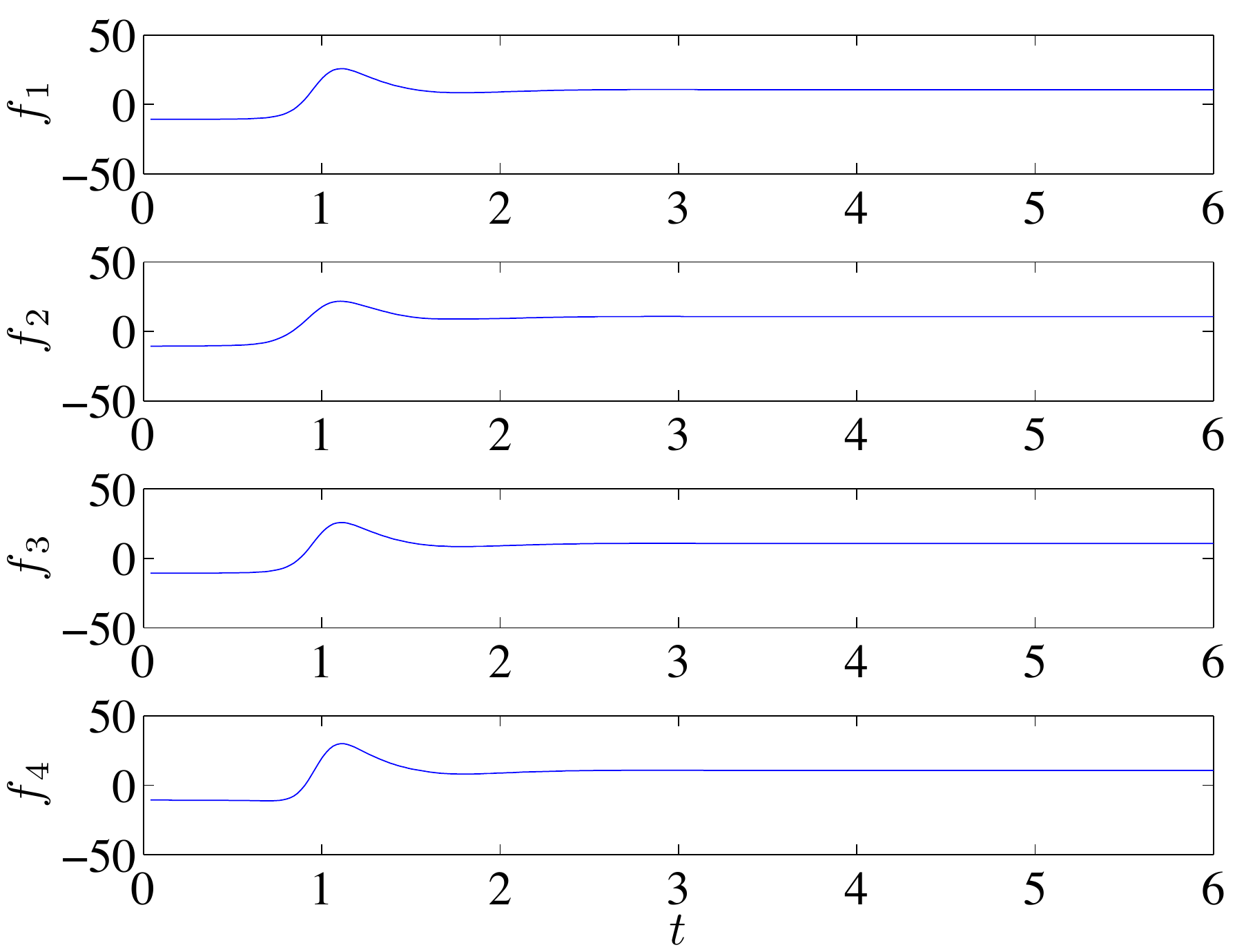}}
}
\caption{Case I: position controlled flight mode for a hovering, recovering from an initially upside down attitude}\label{fig:II}
\end{figure}

Therefore, we cannot apply Proposition \ref{prop:Pos} that gives exponential stability when $\Psi(0)<1$, but by Proposition \ref{prop:Pos2}, we can guarantee exponential attractiveness. From Proposition \ref{prop:Att}, the \textit{attitude error} function $\Psi$ decreases; it eventually becomes less than 1 at $t=0.88$ seconds as illustrated in Figure \ref{fig:IIPsi}. At that instant, the \textit{attitude error} enters the region of attraction specified in Proposition \ref{prop:Pos}. Therefore, for $t>0.88$ seconds, the position tracking error converges to zero exponentially as shown in Figures \ref{fig:IIx}. The region of attraction of the proposed control system almost covers $\SO$, so that the controlled quadrotor UAV can recover from being initially upside down. 

\begin{figure*}
\setlength{\unitlength}{0.70\textwidth}\footnotesize
\centerline{
\begin{picture}(1,0.4)(0,0)
\put(0,0){\includegraphics[width=0.70\textwidth]{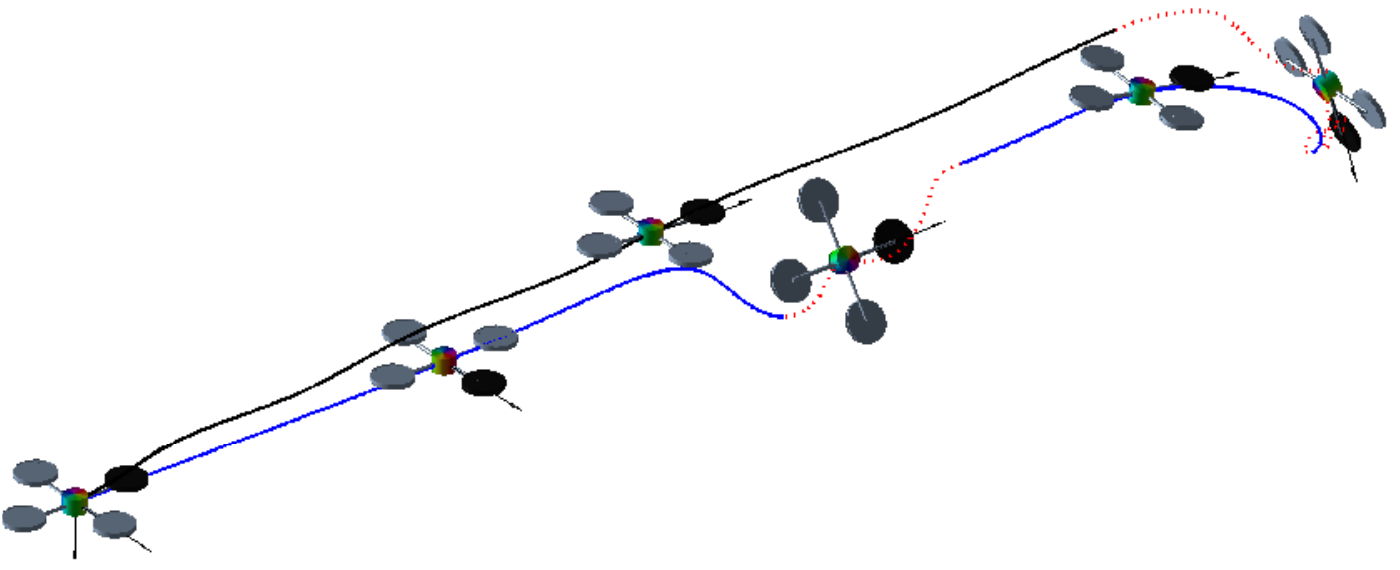}}
\put(-0.01,0.1){\shortstack[c]{initial/terminal\\ position}}
\put(0.36,0.30){\shortstack[c]{(a) velocity tracking}}
\put(0.96,0.36){\shortstack[c]{(b) attitude tracking}}
\put(0.76,0.25){\shortstack[c]{(c) position\\ tracking}}
\put(0.55,0.1){\shortstack[c]{(d) attitude tracking}}
\put(0.28,0.06){\shortstack[c]{(e) position\\ tracking}}
\put(0.11,0.05){\shortstack[c]{$\vec e_1$}}
\put(0.09,-0.008){\shortstack[c]{$\vec e_2$}}
\put(0.025,-0.005){\shortstack[c]{$\vec e_3$}}
\put(0.54,0.25){\shortstack[c]{$\vec b_1$}}
\put(0.87,0.32){\shortstack[c]{$\vec b_1$}}
\put(0.37,0.12){\shortstack[c]{$\vec b_1$}}
\end{picture}}
\caption{Case II: complex maneuver of a quadrotor UAV involving a rotation by $720^\circ$ about $\vec e_2$ (b), and a rotation by $360^\circ$ about $\vec e_1$ (d), with transitions between several flight modes. The direction of the first body-fixed axis is specified for velocity/position tracking modes ((a),(c),(e)) (an animation illustrating this maneuver is available at http://my.fit.edu/\~{}taeyoung).}\label{fig:III3d}
\end{figure*}
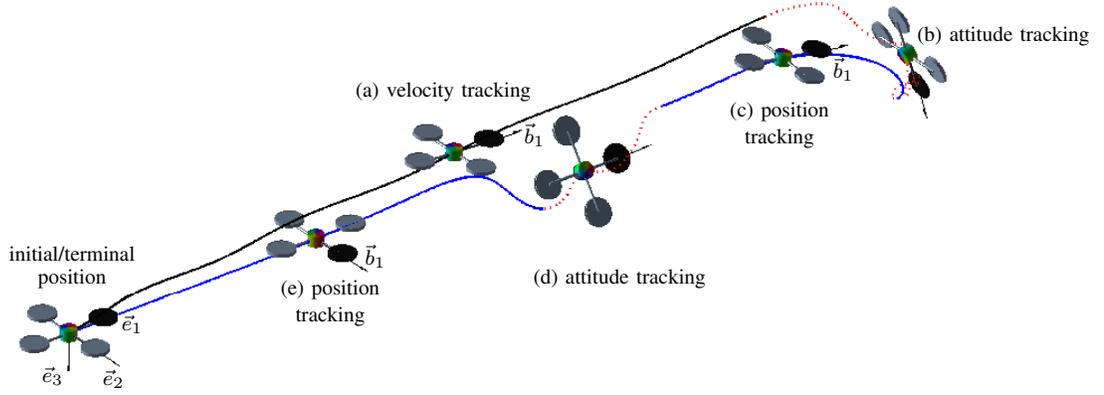

\paragraph*{Case (II): Transition Between Several Flight Modes}
This flight maneuver consists of a sequence of five flight modes, including a rotation by $720^\circ$ (see \reffig{III3d}). 
\begin{itemize}
\item[(a)] Velocity controlled flight mode ($t\in[0,4)$)
\begin{gather*}
v_d(t)=[1+0.5t,\ 0.2\sin(2\pi t),\, -0.1],\\ b_{1_d} (t) = [1,0 ,0].
\end{gather*}
\item[(b)] Attitude controlled flight mode ($t\in[4,6)$): rotation about the second body-fixed axis 
by $720^\circ$
\begin{gather}
R_d(t) = \exp (2\pi(t-4)\hat e_2).\nonumber
\end{gather}
\item[(c)] Position controlled flight mode ($t\in[6,8)$)
\begin{gather*}
x_d(t) = [14-t,\,0,\,0],\quad  b_{1_d} (t) = [1,0 ,0].
\end{gather*}
\item[(d)] Attitude controlled flight mode ($t\in[8,9)$): rotation about the first body-fixed axis
by $360^\circ$
\begin{gather}
R_d(t) = \exp (2\pi(t-8)\hat e_1).\nonumber
\end{gather}
\item[(e)] Position controlled flight mode ($t\in[9,12]$)
\begin{gather*}
x_d(t) = [20-\frac{5}{3}t,\,0,\,0],\quad b_{1_d} (t) = [0, 1, 0].
\end{gather*}
\end{itemize}
Initial conditions are same as the first case.

\begin{figure}
\centerline{
	\subfigure[Attitude error function $\Psi$]{
		\includegraphics[width=0.48\columnwidth]{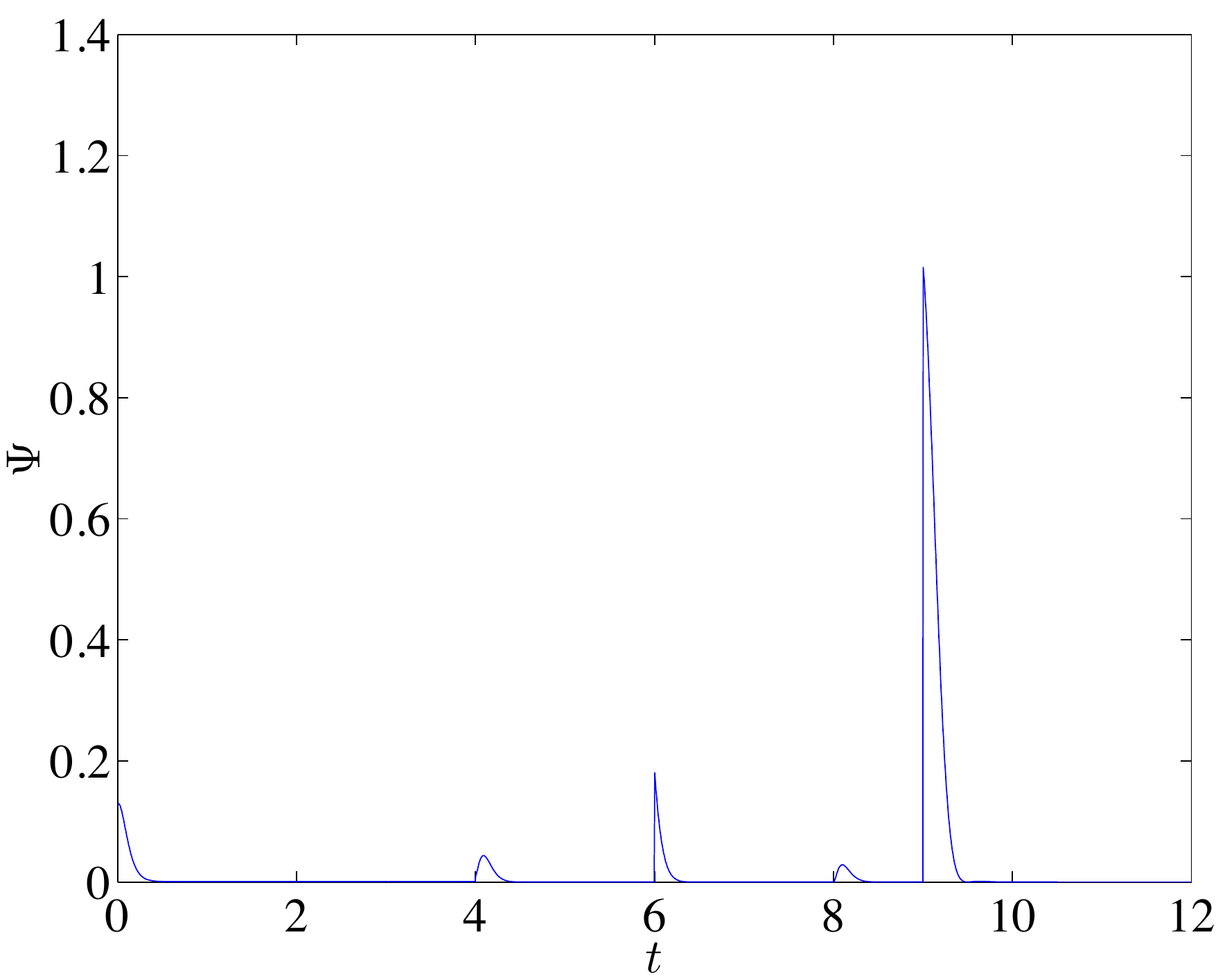}\label{fig:IIIPsi}}
		\hfill
	\subfigure[Position ($x$:solid, $x_d$:dotted, ($\mathrm{m}$))]{
		\includegraphics[width=0.49\columnwidth]{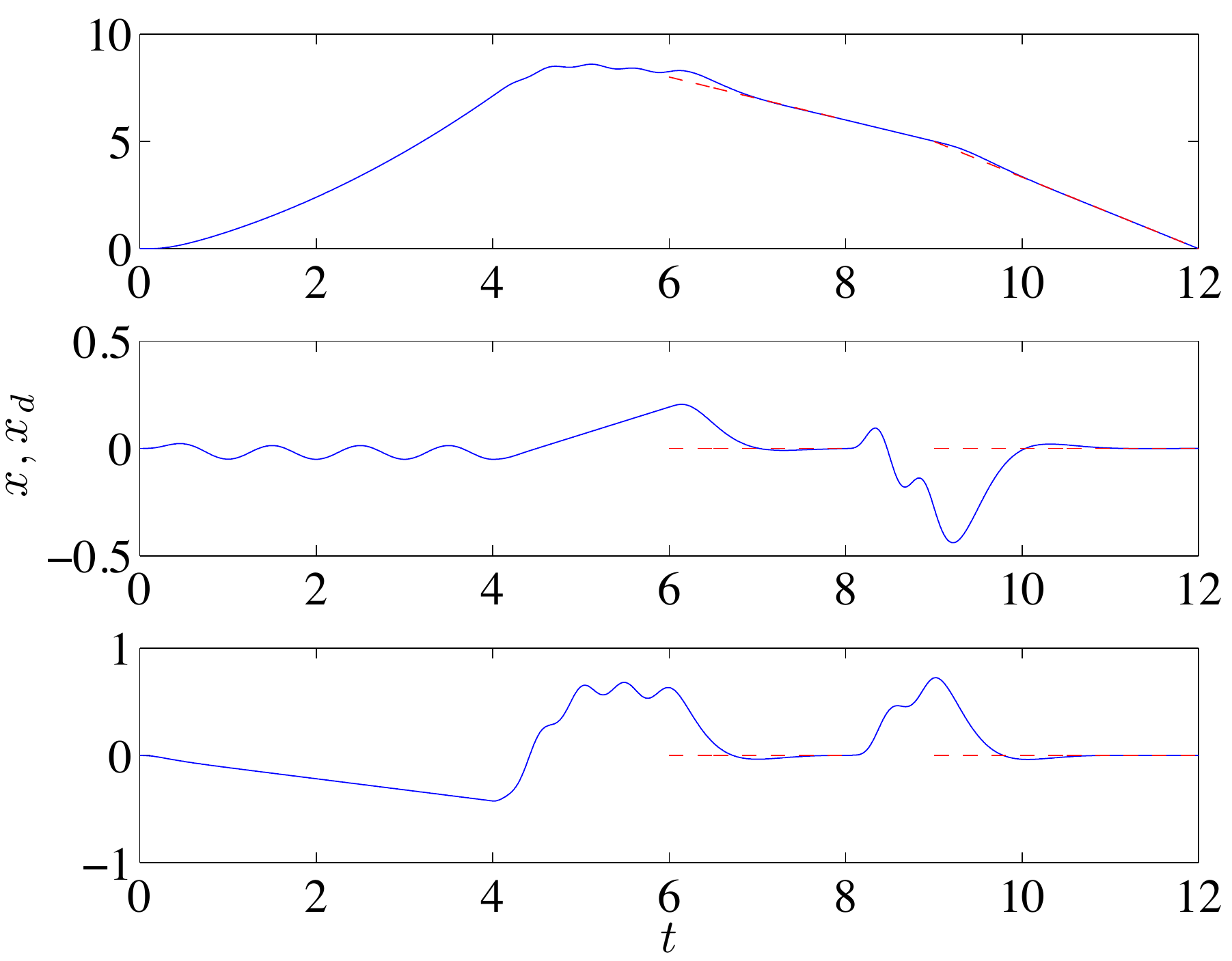}\label{fig:IIIx}}
}
\centerline{
	\subfigure[Angular velocity ($\Omega$:solid, $\Omega_d$:dotted, ($\mathrm{rad/sec}$))]{
		\includegraphics[width=0.49\columnwidth]{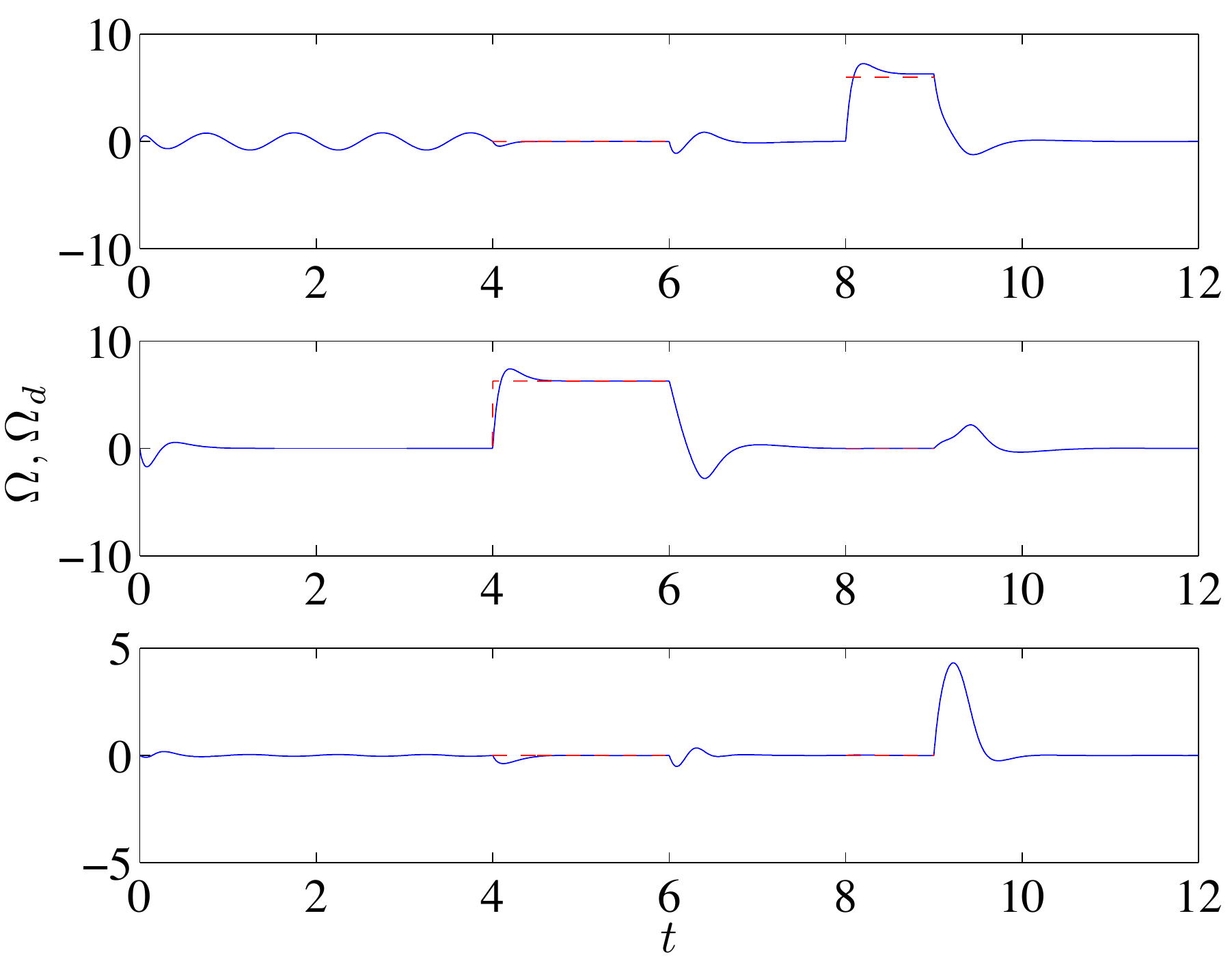}\label{fig:IIIW}}
		\hfill
	\subfigure[Velocity ($v$:solid, $v_d$:dotted, ($\mathrm{m/s}$))]{
		\includegraphics[width=0.48\columnwidth]{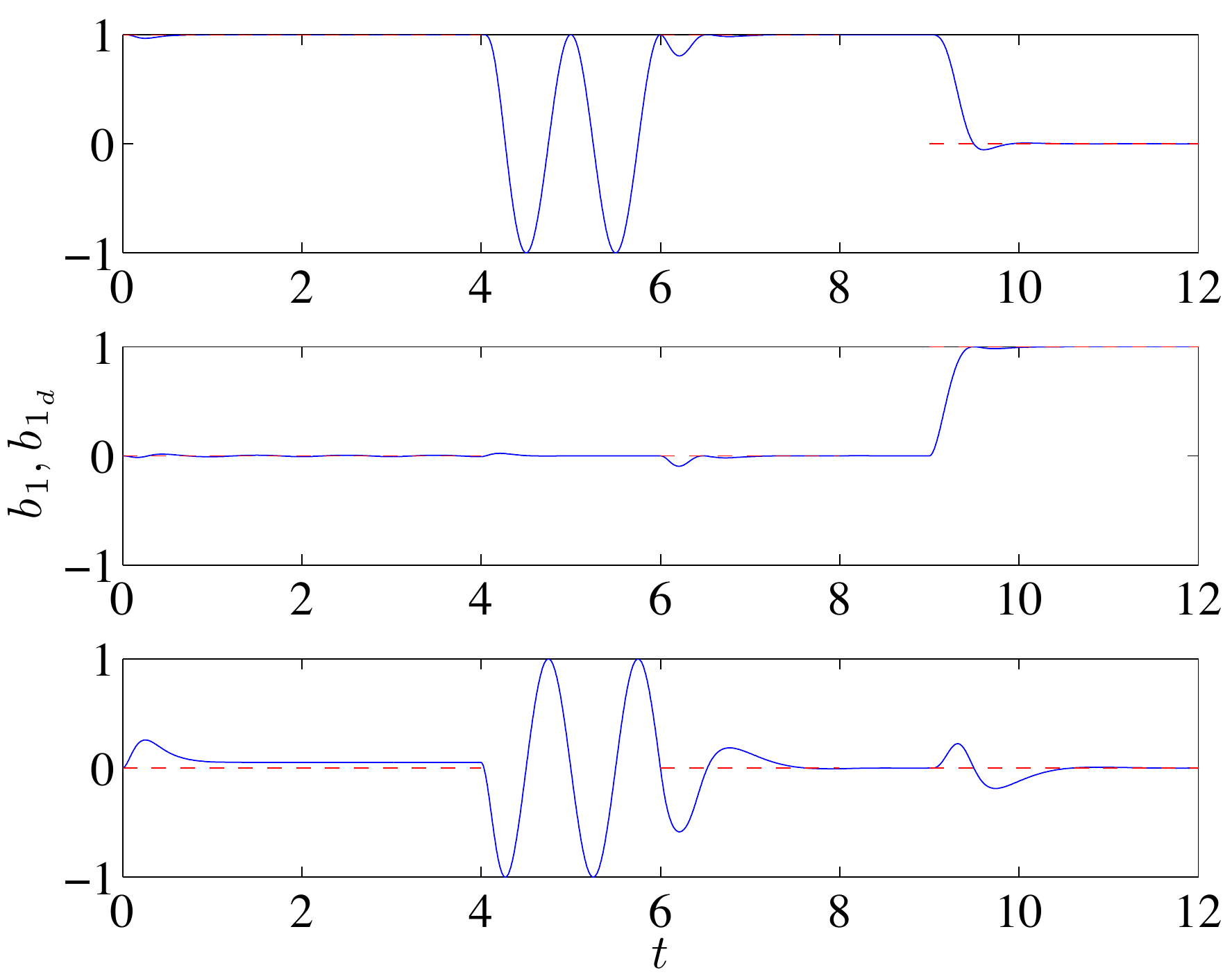}\label{fig:IIIv}}
}
\centerline{
	\subfigure[First body-fixed axis ($b_1$:solid, $b_{1_d}$:dotted)]{
		\includegraphics[width=0.49\columnwidth]{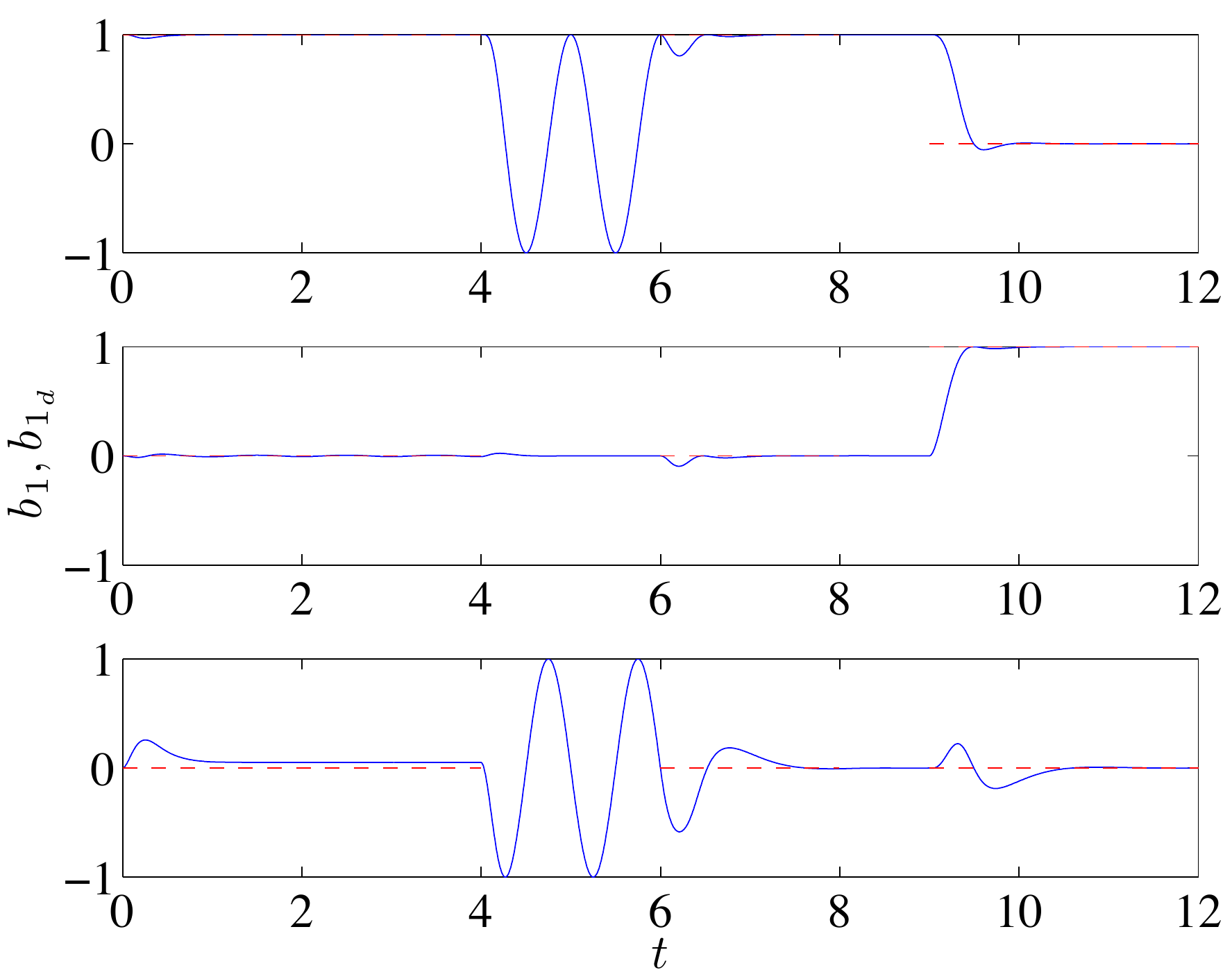}\label{fig:IIIb1}}
}
\caption{Case II: transitions between several flight modes for a complex maneuver}\label{fig:III}
\end{figure}

The second case involves transitions between several flight modes. It begins with a velocity controlled flight mode. As the initial attitude error function is less than 1, the velocity tracking error exponentially converges as shown at \reffig{IIIv}, and the first body-fixed axis asymptotically lies in the plane spanned by $b_{1_d}=e_1$ and $ge_3-\dot v_d$. Since $\|\dot v_d\|\ll g$, the first body-fixed axis remains close to the plane composed of $e_1$ and $e_3$, as illustrated in \reffig{IIIb1}.

This is followed by an attitude tracking mode to rotate the quadrotor by $720^\circ$ about the second body-fixed axis 
according to Proposition \ref{prop:Att}. As discussed in Section \ref{sec:ACFM}, the thrust magnitude $f$ can be arbitrarily chosen in an attitude controlled flight mode.  We cannot apply the results of Proposition \ref{prop:Alt} for altitude tracking, since the third body-fixed axis becomes horizontal several times during the given attitude maneuver. Here we choose the thrust magnitude given by
\begin{gather*}
f(t) = ( k_x (x(t)-x_c) + k_v v(t) + mg e_3)\cdot R(t)e_3,\nonumber
\end{gather*}
which is equivalent to the thrust magnitude for the position controlled flight mode given in \refeqn{f}, when $x_d(t)=x_c=[8,0,0]$. This does not guarantee asymptotic convergence of the quadrotor UAV position to $[8,0,0]$ since the direction of the total thrust is determined by the given attitude command. But, it has the effects that the position of the quadrotor UAV stays close to $x_c$, as illustrated at \reffig{IIIx}. 


Next, a position tracking mode is again engaged, and the quadrotor UAV soon follows a straight line. Another attitude tracking mode and a position tracking mode are repeated to rotate the quadrotor by $360^\circ$ about the direction of the second body-fixed axis.
The thrust magnitude is chosen as 
\begin{gather*}
f(t) = ( k_x (x(t)-x_c) + k_v v(t) + mg e_3)\cdot R(t)e_3,\nonumber
\end{gather*}
where $x_c = [6,0,0]$,
to make the position of the quadrotor UAV remain close to $x_c$ during this attitude maneuver, as discussed above. For the position tracking modes (c) and (e), we have $\ddot x_d=0$, and $b_{1_d}$ lies in the horizontal plane. Therefore, according to Proposition \ref{prop:34C}, the first body-fixed $b_1$ asymptotically converges to $b_{1_d}$, as shown at \reffig{IIIb1}. For example, during the last position tracking mode (e), the first body-fixed axis points to the left of the flight path since $b_{1_d}$ is specified to be $e_2$. These illustrate that by switching between an attitude mode and a position and heading flight mode, the quadrotor UAV can perform the prescribed complex acrobatic maneuver.

\section{CONCLUSIONS}

We presented a global dynamic model for a quadrotor UAV, and we developed tracking controllers for three different flight modes; these were developed in terms of the  special Euclidean group that is intrinsic and coordinate-free, thereby avoiding the singularities of Euler angles and the ambiguities of quaternions in representing attitude.   Using the proposed geometric based controllers for the three flight modes we studied, the quadrotor exhibits exponential stability when the initial attitude error is less than $90^\circ$, and it yields almost global exponentially attractiveness when the initial attitude error is less than $180^\circ$. By switching between different controllers for these flight modes, we have demonstrated that the quadrotor UAV can perform complex acrobatic maneuvers.  Several different complex flight maneuvers were demonstrated in the numerical examples.

\appendix

\section{Properties and Proofs}

\subsection{Properties of the \textit{Hat} Map}\label{app:hat}
\noindent The hat map $\hat\cdot :\Re^3\rightarrow\so$ is defined as
\begin{align}
    \hat x = \begin{bmatrix} 0 & -x_3 & x_2\\
                                x_3 & 0 & -x_1\\
                                -x_2 & x_1 & 0 \end{bmatrix}\label{eqn:hat}
\end{align}
for $x=[x_1;x_2;x_3]\in\Re^3$. This identifies the Lie algebra $\so$ with $\Re^3$ using the vector cross product in $\Re^3$. The inverse of the hat map is referred to as the \textit{vee} map, $\vee:\so\rightarrow\Re^3$. Several properties of the hat map are summarized as follows.
\begin{gather}
    \hat x y = x\times y = - y\times x = - \hat y x,\\
    -\frac{1}{2}\tr{\hat x \hat y} = x^T y,\\
    \tr{\hat x A}=\tr{A\hat x }=\frac{1}{2}\tr{\hat x (A-A^T)}=-x^T (A-A^T)^\vee,\label{eqn:hat1}\\
    \hat x  A+A^T\hat x=(\braces{\tr{A}I_{3\times 3}-A}x)^{\wedge},\label{eqn:xAAx}\\
R\hat x R^T = (Rx)^\wedge,\label{eqn:hat2}
\end{gather}
for any $x,y\in\Re^3$, $A\in\Re^{3\times 3}$, and $R\in\SO$.

\subsection{Proof of Proposition \ref{prop:Att}}\label{sec:pfAtt}

We first find the error dynamics for $e_R,e_\Omega$, and define a Lyapunov function. Then, we show that under the given conditions, $R(t)$ always lies in the sublevel set $L_2$, which guarantees the positive-definiteness of the attitude error function $\Psi$. From this, we show exponential stability of the attitude error dynamics. 

\paragraph{Attitude Error Dynamics}
We find the error dynamics for $\Psi,e_R,e_\Omega$ as follows. Using the attitude kinematics equations, namely $\dot R=R\hat\Omega$, $\dot R_d = R_d\hat\Omega_d$, and equation \refeqn{hat2}, the time derivative of $\Psi$ is given by
\begin{align*}
\dot \Psi(R,R_d) & 
=-\frac{1}{2}\tr{-\hat\Omega_d R_d^T R + R_d^T R\hat\Omega}\\
& = -\frac{1}{2}\tr{R_d^T R(\hat\Omega-R^T R_d\hat\Omega_d R_d^T R)}.
\end{align*}
By \refeqn{eW}, \refeqn{hat1}, this can be written as
\begin{align}
\dot \Psi(R,R_d) & =\frac{1}{2}e_\Omega^T (R_d^T R -R^T R_d)^\vee = e_R\cdot e_\Omega.\label{eqn:Psidot}
\end{align}

Using equations \refeqn{Redot} and \refeqn{xAAx}, the time derivative of $e_R$ can be written as
\begin{align}
\dot e_R & = \frac{1}{2} ( R_d^T R \hat e_\Omega +\hat e_\Omega R^T R_d)^\vee\nonumber\\
& = \frac{1}{2}(\tr{R^T R_d}I-R^T R_d) e_\Omega \equiv C(R_d^TR) e_\Omega.\label{eqn:eRdot}
\end{align}
Now we show that $\|C(R_d^T R)\|_2 \leq 1$ for any $R_d^T R\in\SO$. Using Rodrigues' formula~\cite{BulLew05}, 
we can show that the eigenvalues of $C^T(\exp \hat x) C(\exp \hat x)$ are given by $\cos^2\norm{x}$, $\frac{1}{2}(1+\cos\norm{x})$, and $\frac{1}{2}(1+\cos\norm{x})$, which are less than or equal to $1$ for any $x\in\Re^3$. Therefore, 
 $\|C(R_d^T R)\|_2 \leq 1$, and this implies that 
\begin{align}
\norm{\dot e_R}\leq \norm{e_\Omega}.\label{eqn:neRdot}
\end{align}

From equation \refeqn{eW}, the time derivative of $e_\Omega$ is given by
\begin{align*}
J\dot e_\Omega = J\dot\Omega +J(\hat\Omega R^T R_d \Omega_d - R^T R_d\dot\Omega_d),
\end{align*}
where we use a property of the hat map, $\hat x x =0$ for any $x\in\Re^3$. Substituting the equation of motion \refeqn{EL4} and the control moment \refeqn{M}, this reduces to
\begin{align}
J\dot e_\Omega = -k_R e_R -k_\Omega e_\Omega.\label{eqn:eWdot}
\end{align}
In short, the attitude error dynamics are given by equations \refeqn{Psidot}, \refeqn{eRdot}, \refeqn{eWdot}, and they satisfy \refeqn{neRdot}.

\paragraph{Lyapunov Candidate}
For a non-negative constant $c_2$, let a Lyapunov candidate $\mathcal{V}_2$ be 
\begin{align}
\mathcal{V}_2 = \frac{1}{2} e_\Omega \cdot J e_\Omega + k_R\, \Psi(R,R_d)+c_2 e_R\cdot e_\Omega.\label{eqn:V2}
\end{align}
From equations \refeqn{Psidot}, \refeqn{eRdot}, \refeqn{eWdot}, the time derivative of $\mathcal{V}_2$ is given by
\begin{align}
\dot{\mathcal{V}}_2  & = e_\Omega \cdot J\dot e_\Omega +k_R e_R\cdot e_\Omega
+c_2 \dot e_R \cdot e_\Omega + c_2 e_R \cdot \dot e_\Omega\nonumber\\
& = -k_\Omega \|e_\Omega\|^2 -c_2k_R e_R\cdot J^{-1}e_R +c_2 C(R_d^TR)e_\Omega \cdot e_\Omega\nonumber\\
&\quad  - c_2k_\Omega e_R\cdot J^{-1} e_\Omega.\label{eqn:V2dot}
\end{align}
Since $\|C(R_d^TR)\|\leq 1$, this is bounded by
\begin{align}
\dot{\mathcal{V}}_2   
& \leq - z_2^T W_2 z_2,\label{eqn:V2dot1}
\end{align}
where $z_2 =[\|e_R\|,\;\|e_\Omega\|]^T$, and the matrix $W_2\in\Re^{2\times 2}$ is given by
\begin{align}
W_2 = \begin{bmatrix} \frac{c_2k_R}{\lambda_{M}(J)} & -\frac{c_2k_\Omega}{2\lambda_{m}(J)} \\ 
-\frac{c_2k_\Omega}{2\lambda_{m}(J)} & k_\Omega-c_2 \end{bmatrix}.\label{eqn:W2}
\end{align}

\paragraph{Boundedness of $\Psi$} 
Define $\mathcal{V}'_2 = \mathcal{V}_2\big|_{c_2=0}$. From \refeqn{V2}, \refeqn{V2dot}, we have
\begin{align*}
\mathcal{V}'_2 &  = \frac{1}{2} e_\Omega \cdot J e_\Omega + k_R\, \Psi(R,R_d),\\
\dot{\mathcal{V}}'_2  & = -k_\Omega \|e_\Omega\|^2 \leq 0.
\end{align*}
This implies that $\mathcal{V}'_2$ is non-increasing, i.e., $\mathcal{V}'_2 (t) \leq \mathcal{V}'_2 (0)$.  Using \refeqn{eWb}, the initial value of $\mathcal{V}'_2$ is bounded by $\mathcal{V}'_2 (0) < 2k_R$. Therefore, we obtain
\begin{align}
k_R \Psi(R(t),R_d(t)) \leq \mathcal{V}'_2 (t) \leq \mathcal{V}'_2 (0) < 2 k_R.\label{eqn:kRPsi}
\end{align}
Therefore, the attitude error function is bounded by
\begin{align}
\Psi(R(t),R_d(t)) \leq \psi_2 < 2,\quad \text{for any $t\geq 0$},\label{eqn:Psib0}
\end{align}
and for $\psi_2=\frac{1}{k_R}\mathcal{V}'_2(0)$. Therefore, $R(t)$ always lies in the sublevel set $L_2=\{R\in\SO\,|\,\Psi(R,R_d)<2\}$.

\paragraph{Exponential Stability} Now, we show exponential stability of the attitude dynamics by considering the general case where the constant $c_2$ is positive. Using Rodrigues' formula, we can show that 
\begin{align}
\Psi(R,R_d) &= 1-\cos\norm{x},\label{eqn:Psix}\\
\norm{e_R}^2 & = \sin^2\norm{x}=(1+\cos\norm{x})\Psi(R,R_d)\nonumber\\
& =(2-\Psi(R,R_d))\Psi(R,R_d),\label{eqn:eR2x}
\end{align}
when $R_d^T R = \exp\hat x $ for $x\in\Re^3$. Therefore, from \refeqn{Psib0}, the attitude error function satisfies
\begin{align}
\frac{1}{2} \norm{e_R}^2 \leq  \Psi(R,R_d) \leq \frac{1}{2-\psi_2} \norm{e_R}^2\label{eqn:eRPsi}.
\end{align}
This implies that $\Psi$ is positive-definite and decrescent. It follows that the Lyapunov function $\mathcal{V}_2$ is bounded as
\begin{gather}
z_2^T M_{21} z_2 \leq \mathcal{V}_2 \leq z_2^T M_{22} z_2,
\label{eqn:V2b}
\end{gather}
where
\begin{align}
M_{21} = \frac{1}{2}\begin{bmatrix} k_R & -c_2 \\ -c_2 & \lambda_{m}(J)  \end{bmatrix},\,
M_{22} = \frac{1}{2}\begin{bmatrix} \frac{2k_R}{2-\psi_2} & c_2 \\ c_2 & \lambda_{M}(J)\end{bmatrix}.
\label{eqn:M2}
\end{align}
We choose the positive constant $c_2$ such that
\begin{align*}
c_2 < \min\bigg\{ & k_\Omega,\frac{4k_\Omega k_R\lambda_{m}(J)^2}{k_\Omega^2\lambda_{M}(J)+4k_R\lambda_{m}(J)^2},\sqrt{k_R\lambda_{m}(J)}\bigg\},
\end{align*}
which makes the matrix $W_2$ in \refeqn{V2dot1} and the matrices $M_{21},M_{22}$ in \refeqn{V2b} positive-definite. Therefore, we obtain
\begin{gather}
\lambda_{m}(M_{21})\|z_2\|^2 \leq \mathcal{V}_2 \leq \lambda_{M}(M_{22}) \|z_2\|^2,\\
\dot{\mathcal{V}}_2\leq -\lambda_{m}(W_2) \|z_2\|^2. \label{eqn:V2bb}
\end{gather}
Let $\beta_2=\frac{\lambda_{m}(W_2)}{\lambda_{M}(M_{22})}$. Then, we have
\begin{align}
\dot{\mathcal{V}}_2 \leq - \beta_2\mathcal{V}_2. \label{eqn:V2dotbb}
\end{align}
Therefore, the zero equilibrium of the attitude tracking error $e_R,e_\Omega$ is exponentially stable. Using \refeqn{eRPsi}, this implies that
\begin{align*}
(2-\psi_2)& \lambda_{m}(M_{21}) \Psi\leq  \lambda_{m}(M_{21})\|e_R\|^2 \\
& \leq \lambda_{m}(M_{21})\|z_2\|^2 \leq \mathcal{V}_2(t) \leq \mathcal{V}_2(0) e^{-\beta_2 t}.
\end{align*}
Thus, the attitude error function $\Psi$ exponentially decreases. But, from \refeqn{Psib0}, it is also guaranteed that $\Psi < 2$. This yields \refeqn{Psib}.

\subsection{Proof of Proposition \ref{prop:Alt}}\label{sec:pfAlt}

The rotational dynamics \refeqn{EL3}, \refeqn{EL4} are decoupled from the translational dynamics \refeqn{EL1}, \refeqn{EL2}. As the control moment and assumptions are identical to Proposition \ref{prop:Att}, all of the conclusions of Proposition \ref{prop:Att} hold. 

To show altitude tracking, we take the dot product of \refeqn{EL2} with $e_3$ to obtain
\begin{align*}
m\ddot x_{3} = mg - f e_3\cdot R e_3. 
\end{align*}
Substituting \refeqn{af} into this, we obtain the altitude error dynamics as follows:
\begin{align*}
m\ddot x_{3} = -k_x(x_3-x_{3_d}) -k_v(\dot x_3 - \dot x_{3_d}) + m\ddot x_{3_d}.
\end{align*}
It it clear that this second-order linear system is exponentially stable for positive $k_x,k_v$.

\subsection{Proof of Proposition \ref{prop:Pos}}\label{sec:pfPos}
\setcounter{paragraph}{0}

We first derive the tracking error dynamics. Using a Lyapunov analysis, we show that the velocity tracking error is uniformly bounded, from which we establish exponential stability.

\paragraph{Boundedness of $e_R$}
The assumptions of Proposition \ref{prop:Pos}, namely \refeqn{Psi0}, \refeqn{eWb2} imply satisfaction of the assumptions of Proposition \ref{prop:Att}, \refeqn{eRb0}, \refeqn{eWb}, replacing the notation  $R_d$ by $R_c$.  Therefore, the results of Proposition 1 can be directly applied throughout this proof. From \refeqn{eWb2}, equation \refeqn{kRPsi} can be replaced by
\begin{align}
k_R \Psi(R(t),R_c(t)) \leq \mathcal{V}'_2(0) < k_R\psi_1.\label{eqn:kRPsi2}
\end{align}
Therefore, the \textit{attitude error} function is bounded by
\begin{align}
\Psi(R(t),R_d(t)) \leq \psi_1 < 1,\quad \text{for any $t\geq 0$.}\label{eqn:Psibpf2}
\end{align}
This implies that for the attitude always lies in the sublevel set $L_1=\{R\in\SO\,|\, \Psi(R,R_c)< 1\}$. From \refeqn{Psix}, the \textit{attitude error} is less than $90^\circ$.
Similar to \refeqn{eRPsi}, we can show that 
\begin{align}
\frac{1}{2} \norm{e_R}^2 \leq  \Psi(R,R_c) \leq \frac{1}{2-\psi_1} \norm{e_R}^2\label{eqn:eRPsi1}.
\end{align}

\EditTL{
We also define the following domain $D$
\begin{align}
D=\{(e_x,e_v,R,e_\Omega)\in\Re^3\times\Re^3\times L_1\times\Re^3\,|\, \|e_x\|< e_{x_{\max}}\},\label{eqn:D}
\end{align}
for a fixed constant $e_{x_{\max}}$, restricting the magnitude of the position error. The subsequent Lyapunov analysis is developed in this domain $D$.}

\paragraph{Translational Error Dynamics} The time derivative of the position error is $\dot e_x=e_v$. The time-derivative of the velocity error is given by
\begin{align}
m\dot e_v = m\ddot x -m\ddot x_d = mg e_3 - fRe_3 -m\ddot x_d. \label{eqn:evdot0}
\end{align}
Consider the quantity $e_3^T R_c^T R e_3$, which represents the cosine of the angle between $b_3=Re_3$ and $b_{c_3}=R_ce_3$. Since $1-\Psi(R,R_c)$ represents the cosine of the eigen-axis rotation angle between $R_c$ and $R$, as discussed in \refeqn{Psix}, we have $1 > e_3^T R_c^T R e_3> 1-\Psi(R,R_c)>0$. Therefore, the quantity $\frac{1}{e_3^T R_c^T R e_3}$ is well-defined. To rewrite the error dynamics of $e_v$ in terms of the \textit{attitude error} $e_R$, we add and subtract $\frac{f}{e_3^T R_c^T R e_3}R_c e_3$ to the right hand side of \refeqn{evdot0} to obtain
\begin{align}
m\dot e_v &  = mg e_3 -m\ddot x_d- \frac{f}{e_3^T R_c^T R e_3}R_c e_3 - X,\label{eqn:evdot1}
\end{align}
where $X\in\Re^3$ is defined by
\begin{align}
X=\frac{f}{e_3^T R_c^T R e_3}( (e_3^T R_c^T R e_3)R e_3 -R_ce_3).\label{eqn:X}
\end{align}
Let $A=-k_x e_x - k_v e_v -mg e_3 + m\ddot x_d$.
Then, from \refeqn{f}, \refeqn{Rd3}, we have $f=-A\cdot Re_3$ and ${b}_{3_c}=R_c e_3 = -A/\norm{A}$, i.e. $-A=\|A\| R_c e_3$. By combining these, we obtain $f= (\norm{A}R_c e_3)\cdot R e_3$. Therefore, the third term of the right hand side of \refeqn{evdot1} can be written as
\begin{align*}
- \frac{f}{e_3^T R_c^T R e_3}R_c e_3 & = -\frac{(\norm{A}R_c e_3)\cdot R e_3}{e_3^T R_c^T R e_3}\cdot - \frac{A}{\norm{A}}=A\\
& =-k_x e_x - k_v e_v -mg e_3 + m\ddot x_d.
\end{align*}
Substituting this into \refeqn{evdot1}, the error dynamics of $e_v$ can be written as
\begin{align}
m\dot e_v & =  -k_x e_x - k_v e_v - X.\label{eqn:evdot}
\end{align}

\paragraph{Lyapunov Candidate for Translation Dynamics}
For a positive constant $c_1$, let a Lyapunov candidate $\mathcal{V}_1$ be
\begin{align}
\mathcal{V}_1 = \frac{1}{2}k_x\|e_x\|^2  + \frac{1}{2} m \|e_v\|^2 + c_1 e_x\cdot e_v\label{eqn:V1}.
\end{align}
The derivative of ${\mathcal{V}}_1$ along the solution of \refeqn{evdot} is given by
\begin{align}
\dot{\mathcal{V}}_1 & = k_x e_x\cdot e_v  + e_v \cdot \{-k_x e_x -k_v e_v+X\}
+ c_1 e_v\cdot e_v\nonumber\\
&\quad + \frac{c_1}{m} e_x \cdot \{-k_x e_x -k_v e_v+X\}\nonumber\\
& =  -(k_v-c_1) \|e_v\|^2 
- \frac{c_1 k_x}{m} \|e_x\|^2 
-\frac{c_1 k_v}{m} e_x\cdot e_v\nonumber\\
&\quad+X\cdot \braces{ \frac{c_1}{m} e_x + e_v}.\label{eqn:V1dot0}
\end{align}
We find a bound on $X$ using \refeqn{X} as follows. Since $f=\|A\| (e_3^T R_c^T R e_3)$, we have
\begin{align*}
\norm{X} & \leq \|A\|\,\| (e_3^T R_c^T R e_3)R e_3 -R_ce_3\|\\
& \leq( k_x \|e_x\| + k_v \|e_v\| + B)\, \| (e_3^T R_c^T R e_3)R e_3 -R_ce_3\|.
\end{align*}
The last term $\| (e_3^T R_c^T R e_3)R e_3 -R_ce_3\|$ represents the sine of the angle between $b_3=Re_3$ and $b_{c_3}=R_c e_3$, since
\begin{align*}
(b_{3_c}\cdot b_3)b_3 - b_{3_c} = b_{3}\times (b_3\times b_{3_c}).
\end{align*}
From \refeqn{eR2x}, $\|e_R\|$ represents the sine of the eigen-axis rotation angle between $R_c$ and $R$. Therefore, we have $\| (e_3^T R_c^T R e_3)R e_3 -R_ce_3\| \leq \| e_R\|$. From \refeqn{eR2x}, \refeqn{Psibpf2}, it follows that
\begin{align*}
\| (e_3^T R_d^T R e_3)R e_3 -R_de_3\| &\leq \| e_R\| = \sqrt{\Psi(2-\Psi)}\\
& \leq \sqrt{\psi_1 (2-\psi_1)}\equiv\alpha  <1.
\end{align*}
Therefore, $X$ is bounded by
\begin{align}
\norm{X} 
&\leq ( k_x \|e_x\| + k_v \|e_v\| + B) \|e_R\| \nonumber\\
&\leq ( k_x \|e_x\| + k_v \|e_v\| + B) \alpha.\label{eqn:XB}
\end{align}
Substituting this into \refeqn{V1dot0}, 
\begin{align}
\dot{\mathcal{V}}_1 & \leq   -(k_v-c_1) \|e_v\|^2 
- \frac{c_1 k_x}{m} \|e_x\|^2 
-\frac{c_1k_v}{m}e_x\cdot e_v\nonumber\\
& \quad +( k_x \|e_x\| + k_v \|e_v\| + B) \|e_R\| \braces{ \frac{c_1}{m} \|e_x\| + \|e_v\|}\nonumber\\
& \leq   -(k_v(1-\alpha)-c_1) \|e_v\|^2 
- \frac{c_1 k_x}{m}(1-\alpha) \|e_x\|^2 \nonumber\\
&\quad + \frac{c_1k_v}{m}(1+\alpha) \|e_x\|\|e_v\|\nonumber\\
&\quad +  \|e_R\| \braces{\frac{c_1}{m}B \|e_x\| + B\|e_v\|+k_x\|e_x\|\|e_v\|}.\label{eqn:V1dot1}
\end{align}

\EditTL{In the above expression for $\dot{\mathcal{V}}_1$, there is a third-order error term, namely $k_x\|e_R\|\|e_x\|\|e_v\|$. It is possible to choose its upper bound as $k_x\alpha\|e_x\|\|e_v\|$ similar to other terms, but the corresponding stability analysis becomes complicated, and the initial attitude error should be reduced further. Instead, we restrict our analysis to the domain $D$ defined at \refeqn{D}, and an upper bound is chosen as $k_xe_{x_{\max}}\|e_R\|\|e_v\|$.}

%
%

\paragraph{Lyapunov Candidate for the Complete System:}
Let $\mathcal{V}=\mathcal{V}_1+\mathcal{V}_2$ be the Lyapunov candidate of the complete system.
\begin{align}
\mathcal{V} & = \frac{1}{2} k_x \|e_x\|^2 + \frac{1}{2}m \|e_v\|^2 + c_1 e_x\cdot e_v\nonumber\\
&\quad + \frac{1}{2}e_\Omega \cdot Je_  \Omega + k_R\Psi(R,R_d) + c_2 e_R\cdot e_\Omega.\label{eqn:V}
\end{align}
Using \refeqn{eRPsi1}, the bound of the  Lyapunov candidate $\mathcal{V}$ can be written as
\begin{align}
z_1^T M_{11} z_1 + z_2^T M_{21} z_2 \leq \mathcal{V} \leq z_1^T M_{12} z_1 + z_2^T M'_{22} z_2,\label{eqn:Vb}
\end{align}
where $z_1=[\|e_x\|,\;\|e_v\|]^T$, $z_2=[\|e_R\|,\;\|e_\Omega\|]^T\in\Re^2$, and the matrices $M_{11},M_{12},M_{21},M_{22}$ are given by
\begin{gather*}
M_{11} = \frac{1}{2}\begin{bmatrix} k_x & -c_1 \\ -c_1 & m\end{bmatrix},\quad
M_{12} = \frac{1}{2}\begin{bmatrix} k_x & c_1 \\ c_1 & m\end{bmatrix},\\
M_{21} = \frac{1}{2}\begin{bmatrix} k_R & -c_2 \\ -c_2 & \lambda_{m}(J)  \end{bmatrix},\quad
M'_{22} = \frac{1}{2}\begin{bmatrix} \frac{2k_R}{2-\psi_1} & c_2 \\ c_2 & \lambda_{M}(J)\end{bmatrix}.
\end{gather*}

Using \refeqn{V2dot1} and \refeqn{V1dot1}, the time-derivative of $\mathcal{V}$ is given by
\begin{align}
\dot{\mathcal{V}} \leq -z_1^T W_1 z_1  + z_1^T W_{12} z_2 - z_2^T W_2 z_2, \label{eqn:Vdotb}
\end{align}
where $W_1,W_{12},W_2\in\Re^{2\times 2}$ are defined as follows:
\begin{align}
W_1 &= \begin{bmatrix} \frac{c_1k_x}{m}(1-\alpha) & -\frac{c_1k_v}{2m}(1+\alpha)\\
-\frac{c_1k_v}{2m}(1+\alpha) & k_v(1-\alpha)-c_1\end{bmatrix},\\
W_{12}&=\begin{bmatrix}
\frac{c_1}{m}B & 0 \\ B+k_xe_{x_{\max}} & 0\end{bmatrix},\\
W_2 &= \begin{bmatrix} \frac{c_2k_R}{\lambda_{M}(J)} & -\frac{c_2k_\Omega}{2\lambda_{m}(J)} \\ 
-\frac{c_2k_\Omega}{2\lambda_{m}(J)} & k_\Omega-c_2 \end{bmatrix}.
\end{align}

\paragraph{Exponential Stability} 
Under the given conditions \refeqn{c1b}, \refeqn{c2b} of the proposition, all of the matrices $M_{11}$, $M_{12}$, $W_{1}$, $M_{21}$, $M_{22}$, $W_2$, and the Lyapunov candidate $\mathcal{V}$ become positive-definite, and
\begin{align*}
\dot{\mathcal{V}} &\leq -\lambda_{m}(W_1)\|z_1\|^2 +\|W_{12}\|_2 \|z_1\|\|z_2\| - \lambda_{m} (W_{2})\|z_2\|^2.
\end{align*}
The condition given by \refeqn{kRkWb} guarantees that $\dot{\mathcal{V}}$ becomes negative-definite. Therefore, the zero equilibrium of the tracking errors of the complete dynamics is exponentially stable. \EditTL{A (conservative) region of attraction is characterized by a sub-level set of $\mathcal{V}$ contained in the domain $D$, as written at \refeqn{RAz}, as well as \refeqn{eWb2} required for the boundedness of $e_R$.}

\subsection{Proof of Proposition \ref{prop:Pos2}}\label{sec:pfPos2}
The given assumptions \refeqn{eRb3}, \refeqn{eWb3} satisfy the assumption of Proposition \ref{prop:Att}, from which the tracking error $z_2=[\|e_R\|,\|e_\Omega\|]$ is guaranteed to exponentially decreases, and to enter the region of attraction of Proposition \ref{prop:Pos}, given by \refeqn{Psi0}, \refeqn{eWb2}, in a finite time $t^*$. 

Therefore, if we show that the tracking error $z_1=[\|e_x\|,\|e_v\|]$ is bounded in $t\in[0, t^*]$, then the tracking error $z=[z_1,z_2]$ is uniformly bounded for any $t>0$, and it exponentially decreases for $t>t^*$. This yields exponential attractiveness.

The boundedness of $z_1$ is shown as follows. The error dynamics or $e_v$ can be written as
\begin{align*}
m\dot e_v = mg e_3 -fRe_3 -m\ddot x_d.
\end{align*}
Let ${\mathcal{V}}_3$ be a positive-definite function of $\|e_x\|$ and $\|e_v\|$:
\begin{align*}
{\mathcal{V}}_3 = \frac{1}{2}\|e_x\|^2  + \frac{1}{2}m \|e_v\|^2.
\end{align*}
Then, we have $\|e_x\|\leq \sqrt{2\mathcal{V}_3}$, $\|e_v\|\leq \sqrt{\frac{2}{m}\mathcal{V}_3}$. The time-derivative of $\mathcal{V}_3$ is given by
\begin{align*}
\dot{\mathcal{V}}_3 & =  e_x\cdot e_v + e_v\cdot (mg e_3 -f R e_3- m\ddot x_d) \\
& \leq \|e_x\|\|e_v\| + \|e_v\| \|mg e_3 -m\ddot x_d\| + \|e_v\| \|R e_3\| |f|.
\end{align*}
Using \refeqn{B}, \refeqn{f}, we obtain
\begin{align*}
\dot{\mathcal{V}}_3 & \leq \|e_x\|\|e_v\| + \|e_v\| B + \|e_v\| (k_x\|e_x\|+k_v\|e_v\| + B)\\
& = k_v \|e_v\|^2 + (2B+(k_x+1)\|e_x\|)\|e_v\|\\
& \leq d_1 \mathcal{V}_3 + d_2 \sqrt{\mathcal{V}_3},
\end{align*}
where $d_1 = k_v \frac{2}{m} +2(k_x+1)\frac{1}{\sqrt{m}}$, $d_2 = 2B \sqrt{\frac{2}{m}}$. Suppose that $\mathcal{V}_3 \geq 1$ for a time interval $[t_a,t_b]\subset [0,t^*]$. In this time interval, we have $\sqrt{\mathcal{V}_3} \leq \mathcal{V}_3$. Therefore, 
\begin{align*}
\dot{\mathcal{V}}_3 \leq (d_1+d_2) \mathcal{V}_3 \quad \Rightarrow\quad \mathcal{V}_3(t) \leq \mathcal{V}_3(t_a) e^{(d_1+d_2)(t-t_a)}.
\end{align*}
Therefore, for any time interval in which $\mathcal{V}_3\geq 1$, $\mathcal{V}_3$ is bounded. This implies that $\mathcal{V}_3$ is bounded for $0\leq t\leq t^*$.

In summary, for any initial condition satisfying \refeqn{eRb3},\refeqn{eWb3}, the tracking error converges to the region of attraction for exponential stability according to Proposition \ref{prop:Pos}, and during that time period, tracking errors are bounded. Therefore, the zero equilibrium of the tracking error is exponentially attractive.


\end{document}